\documentclass[a4paper,reqno]{amsart}
\usepackage{amsfonts}

\catcode`@=11 
\def\blfootnote{\xdef\@thefnmark{}\@footnotetext}
\catcode`@=12 

\theoremstyle{plain}
\begingroup
\newtheorem{theorem}{Theorem}[section]
\newtheorem{lemma}[theorem]{Lemma}
\newtheorem{proposition}[theorem]{Proposition}

\endgroup

\theoremstyle{definition}
\begingroup
\newtheorem{definition}[theorem]{Definition}
\newtheorem{remark}[theorem]{Remark}

\endgroup

\theoremstyle{remark}
\begingroup

\endgroup

\mathchardef\emptyset="001F

\numberwithin{equation}{section}

\newcommand{\e}{\varepsilon}
\newcommand{\Om}{\Omega}

\newcommand{\R}{{\mathbb R}}

\newcommand{\Z}{\mathbb Z}

\newcommand{\wto}{\rightharpoonup}

\newcommand{\setmeno}{\!\setminus\!}

\renewcommand{\hom}{{hom}}
\newcommand{\supp}{{\rm supp}}

\newcommand{\Lone}{{\mathcal L}^1}

\newcommand{\bary}{{\rm bar}}

\newcommand{\ol}{\overline}

\title[Time-dependent systems of generalized Young measures]
{Time-dependent systems of generalized Young measures}
\author{G.\ Dal Maso}
\author{A.\ DeSimone}
\author{M.G.\ Mora}
\author{M.\ Morini}
\address[G.~Dal Maso, A.~DeSimone, M.G.~Mora, and M.\ Morini]{SISSA, Via Beirut 
4, 
34014 Trieste, Italy}
\email[Gianni Dal Maso]{dalmaso@sissa.it}
\email[Antonio DeSimone]{desimone@sissa.it}
\email[Maria Giovanna Mora]{mora@sissa.it}
\email[Massimiliano Morini]{morini@sissa.it}

\begin{document}
\begin{abstract}
In this paper some new tools for the study of evolution problems in the framework of 
Young measures are introduced. A suitable notion of time-dependent system of 
generalized Young measures is defined, which allows to extend the classical notions of 
total variation and absolute continuity with respect to time, as well as the notion of time 
derivative. The main results are a Helly type theorem for sequences of systems of 
generalized Young measures and a theorem about the existence of the time derivative for 
systems with bounded variation  with respect to time. 
\end{abstract}
\maketitle

{\small

\bigskip
\keywords{\noindent {\bf Keywords:} Young measures, bounded variation, absolute 
continuity, weak derivatives, concentration and oscillation effects \blfootnote{Preprint SISSA 98/2005/M (December 2005)}
}

\subjclass{\noindent {\bf 2000 Mathematics Subject Classification:} 28A33%
, 26A45%
}
}
\tableofcontents
\bigskip
\bigskip

\begin{section}{Introduction}

\par
The notion of Young measure was introduced by L.C. Young in \cite{You1} to describe 
generalized solutions to minimum problems in the calculus of variations. Since then it has 
been applied to several problems  in the calculus of variations, in control theory, in partial 
differential equations, and in mathematical economics. For the general theory of Young 
measures we refer to \cite{Bald}, \cite{Ball}, \cite{Cas-Ray-Val}, \cite[Chapters 2 
and~3]{Gam}, \cite{Kin-Ped}, \cite{Ped}, \cite{Val}, \cite[Chapter~IV]{War}, and 
\cite{You2}. Several applications are devoted to evolution problems (see, e.g., \cite{Dem}, 
\cite{Dip-Maj}, \cite{Mie}, \cite{Mie2}, \cite{Rie}, and~\cite{Tar}).

\par
In this paper we introduce some new tools in the theory of Young measures for the study 
of rate independent evolution problems. To describe the content of this paper, let us 
consider a problem defined on a time interval $I$, with space variable $x$ in a compact 
metric space $X$, and state variable $u$ in a finite dimensional Hilbert space $\Xi$. We 
assume that $X$ is endowed with a given nonnegative Radon measure $\lambda$ with 
$\supp\,\lambda=X$. Given a sequence $u_k=u_k(t,x)$ of functions from $I{\times}X$ 
to $\Xi$, satisfying suitable estimates, it is often possible to extract a subsequence 
converging, for every $t\in I$, to a Young measure $\mu_t$, which encodes information 
on the statistics of the space oscillations of $u_k(t,x)$ at time~$t$.

\par
To simplify the notation, the Young measure $\mu_t$ will always be regarded as a 
measure on $X{\times}\Xi$, whose projection on $X$ coincides with~$\lambda$. In this 
introduction we will never consider the standard disintegration $(\mu_t^x)_{x\in X}$, 
which is usual in the classical presentation of the theory (see 
Remark~\ref{Disintegration}).

\par
If we want to extend some natural notions, like total variation, absolute continuity, or time 
derivative, from the original context of time dependent functions to the generalized context 
of time-dependent Young measures, we need to know the joint oscillations of 
$u_k(t_1,x),\dots,u_k(t_m,x)$ for every finite sequence $t_1,\dots,t_m$ of times. These 
are described by the Young measure $\mu_{t_1\dots t_m}$, with state space $\Xi^m$, 
generated by the sequence of $\Xi^m$-valued functions  $(u_k(t_1,x),\dots,u_k(t_m,x))$. 
It is easy to see that $\mu_{t_1\dots t_m}$ cannot be derived from the measures 
$\mu_{t_1},\dots,\mu_{t_m}$. Indeed, these measures give no information on the 
correlation between the oscillations at different times. The situation is similar to what 
happens in stochastic processes, where the knowledge of the distribution function of each 
single random variable is not enough to deduce their joint distribution.

\par
This leads to the notion of system of Young measures, defined as a family 
$(\mu_{t_1\dots t_m})$, where $t_1,\dots, t_m$ run over all finite sequences of elements 
of $I$, with $t_1<\dots<t_m$, and each $\mu_{t_1\dots t_m}$ is a Young measure on 
$X$ with values in $\Xi^m$. We assume that $(\mu_{t_1\dots t_m})$ satisfies the 
following compatibility condition, which is always satisfied when $\mu_{t_1\dots t_m}$ 
is generated by a sequence of time-dependent functions: if $\{s_1,\dots, s_n\}\subset 
\{t_1,\dots, t_m\}$ and $s_1<\dots<s_n$, then $\mu_{s_1\dots s_n}$ coincides with the 
corresponding projection of $\mu_{t_1\dots t_m}$.

\par
The notions of total variation (Definition~\ref{variation}), absolute continuity 
(Definition~\ref{abscont}), or time derivative (Definition~\ref{weak*der}) can be easily 
defined in the framework of systems of Young measures in such a way that they coincide 
with the standard notions in the case of time-dependent functions. The main result of the 
paper is a version of Helly's Theorem for systems of Young measures 
(Theorem~\ref{Helly}): if $(\mu^k_{t_1\dots t_m})$ has uniformly bounded variation, 
then there exist a system $(\mu_{t_1\dots t_m})$ with bounded variation, a set 
$\Theta\subset I$, with $I\setminus\Theta$ at most countable, and a subsequence, still 
denoted $(\mu^k_{t_1\dots t_m})$, such that $\mu^k_{t_1\dots t_m}\wto \mu_{t_1\dots 
t_m}$ weakly$^*$ for every finite sequence $t_1,\dots, t_m\in\Theta$ with 
$t_1<\dots<t_m$. 

\par
Another important result provides the existence of the time derivative $\dot\mu_t$ for 
almost every $t$ whenever the family  $(\mu_{t_1\dots t_m})$ has bounded variation 
(Theorem~\ref{BV-weakder}).
The variation can be expressed by an integral involving the time derivatives when 
$(\mu_{t_1\dots t_m})$ is absolutely continuous (Theorem~\ref{weakder}). 

\par
Our motivation for studying systems of Young measures stems from the analysis of
quasistatic evolution problems  with nonconvex energies, which arise in the study of 
plasticity with softening  \cite{DM-DeS-Mor-Mor}.
Since in these applications the energy functionals have linear growth in some directions, 
we have to consider the case where the generating sequence $(u_k(t,x))$ is bounded in 
$L^r_\lambda(X;\Xi)$ only for $r=1$. It is well known that in this case Young measures 
should be replaced by more general objects, which take into account concentrations at 
infinity (see  \cite{Dip-Maj}). In  \cite{Ali-Bou} and \cite{Fon-Mue-Ped} this is done by 
considering a pair $(\mu^Y\!,\mu^\infty)$, where $\mu^Y\!$ is a Young measure on $X$ 
with values in $\Xi$ and $\mu^\infty$, called the varifold measure, is a measure supported 
on $X{\times}\Sigma_{\Xi}$, where $\Sigma_{\Xi}$ denotes the unit sphere in $\Xi$. 

\par
In the spirit of \cite{Dip-Maj}, we prefer to present these generalized Young measures in a 
different way, using homogeneous coordinates to describe the completion of $\Xi$ 
obtained by adding a point at infinity for each direction. We replace the pair  
$(\mu^Y\!,\mu^\infty)$ by a single nonnegative measure $\mu$ on 
$X{\times}\Xi{\times}\R$ (Definition~\ref{GY}), acting only on continuous functions 
$f(x,\xi,\eta)$ which are positively  homogeneous of degree one in $(\xi,\eta)$. We 
assume that $\mu$ is supported on the set $\{\eta\ge 0\}$ and that the projection of 
$\eta\mu$ onto $X$ coincides with $\lambda$. We show that, if $\lambda$ is nonatomic, 
then the space $L^1_\lambda(X;\Xi)$ can be identified (Definition~\ref{GYMu}) with a 
dense subset of the space of generalized Young measures (Theorem~\ref{thm:density}).

\par
Using this approach, we are able to prove the results on total variation and time derivatives  
for systems of Young measures in a context that is general enough for the applications 
considered in \cite{DM-DeS-Mor-Mor}.

\par
\end{section}

\begin{section}{A space of homogeneous functions and its dual}

If $E$ is a locally compact space with a countable base and $\Xi$ is a finite dimensional 
Hilbert space, $M_b(E;\Xi)$ denotes the space of bounded Radon measures on $E$ with 
values in $\Xi$, endowed with the norm $\|\nu\| :=|\nu|(E)$, where $|\nu|$ denotes the 
variation of $\nu$. When $\Xi=\R$, the corresponding space will be denoted simply by 
$M_b(E)$. As usual, $M_b^+(E)$ denotes the cone of nonnegative bounded Radon 
measures on $E$. If $\nu\in M_b(E)$ and $f\in L^1_\nu(E;\Xi)$, the measure $f\,\nu\in 
M_b(E;\Xi)$ is defined by $(f\,\nu)(A):=\int_A f\,d\nu$ for every Borel set $A\subset E$.

\par
By the Riesz Representation Theorem $M_b(E;\Xi)$ can be identified with the dual of 
$C_0(E;\Xi)\!$, the space of continuous functions $\varphi\colon E\to\Xi$ such that 
$\{|\varphi|\ge\e\}$ is compact for every $\e>0$. The weak$^*$ topology of 
$M_b(E;\Xi)$ is defined using this duality.

\par
Throughout the paper $(X,d)$ is a given compact metric space and $\lambda$ is a fixed 
nonnegative Radon measure on $X$ with $\supp\,\lambda=X$. The symbol $\Xi$ will 
denote any finite dimensional Hilbert space. The spaces $L^r(X;\Xi)$, $r\ge 1$, will 
always refer to the measure $\lambda$. If $\mu\in M_b(X;\Xi)$, $\mu^a$ and $\mu^s$ 
denote the absolutely continuous and the singular part of $\mu$ with respect to 
$\lambda$.
Measures in $M_b(X;\Xi)$ which are absolutely continuous with respect to $\lambda$ 
will always be identified with their densities, which belong to $L^1(X;\Xi)$. In this way 
$L^1(X;\Xi)$ is regarded as a subspace of $M_b(X;\Xi)$.

\par
In order to define the notion of generalized Young measure on $X$ with values in $\Xi$, 
it is convenient to introduce a space of homogeneous functions and to discuss some 
properties of its dual.

\par
\begin{definition}\label{Chom}
Let $C^\hom(X{\times}\Xi)$ be the space of all $f\colon X{\times}\Xi\to\R$ such that  
$\xi\mapsto f(x,\xi)$ is positively homogeneous of degree one on $\Xi$ for every $x\in 
X$; i.e., $f(x,t\xi)=tf(x,\xi)$ for every $x\in X$, $\xi\in\Xi$, and $t\ge 0$. This space is 
endowed with the norm
$$
\|f\|_{\hom}:= \max\{ |f(x,\xi)|: x\in X,\ \xi\in\Sigma_\Xi\} \,,
$$
where $\Sigma_\Xi:=\{\xi\in\Xi: |\xi|=1\}$. 
\end{definition}

\par
We introduce now two dense subspaces of $C^\hom(X{\times}\Xi)$ that will be useful in 
the proof of some properties of generalized Young measures.

\par
\begin{definition}\label{ChomL}
Let $C^\hom_L(X{\times}\Xi)$ be the space of all $f\in C^\hom(X{\times}\Xi)$ which 
satisfies the following Lipschitz condition: there exists a constant $a\in\R$ such that
\begin{equation}\label{Lip}
|f(x,\xi_1)-f(x,\xi_2)|\le a\, |\xi_1-\xi_2|
\end{equation}
for every $x\in X$ and every $\xi_1,\xi_2\in\Xi$.
\end{definition}

\par
\begin{remark}\label{rem:ChomL}
If $f\in C^\hom_L(X{\times}\Xi)$ and $\omega$ is the modulus of continuity of the 
restriction of $f$ to $X{\times}\Sigma_\Xi$, then \eqref{Lip} and the homogeneity of 
$f$ imply that
\begin{eqnarray*}
&|f(x_1,\xi_1)-f(x_2,\xi_2)|\le |f(x_1,\xi_1)-f(x_1,\xi_2)| + |f(x_1,\xi_2)-f(x_2,\xi_2)|\le
\\
&\le a\, |\xi_1-\xi_2| + |\xi_2|\,\omega(d(x_1,x_2))\,.
\end{eqnarray*}
Exchanging the roles of $\xi_1$ and $\xi_2$ we obtain
$$
|f(x_1,\xi_1)-f(x_2,\xi_2)|\le a\, |\xi_1-\xi_2| + \min\{|\xi_1|,|\xi_2|\}\,\omega(d(x_1,x_2))
$$
for every $x_1,x_2\in X$ and every $\xi_1,\xi_2\in\Xi$.
\end{remark}

\par
\begin{lemma}\label{L-density}
The space $C^\hom_L(X{\times}\Xi)$ is dense in $C^\hom(X{\times}\Xi)$.
\end{lemma}

\par
\begin{proof}
Let us fix $f\in C^\hom(X{\times}\Xi)$. For every $k>\|f\|_\hom$ let us consider the 
Moreau-Yosida approximation $f_k\colon X{\times}\Xi\to\R$ defined by 
$$
f_k(x,\xi):=\min_{\xi'\in \Xi}\{ f(x,\xi')+k\, |\xi'-\xi| \}\,.
$$
Using the standard properties of Moreau-Yosida approximations it is easy to check that 
$f_k\in C^\hom_L(X{\times}\Xi)$ (with constant $k$) and that the sequence $f_k$ is 
nondecreasing and converges pointwise to $f$ (see, e.g., \cite[Remark~9.6 and 
Theorem~9.13]{DM}). By Dini's Theorem we conclude that $f_k\to f$ uniformly on 
$X{\times}\Sigma_\Xi$, hence $f_k\to f$ in $C^\hom(X{\times}\Xi)$.
\end{proof}

\par
\begin{definition}\label{ChomT}
Let $C^\hom_\triangle(X{\times}\Xi)$ be the space of all $f\in C^\hom(X{\times}\Xi)$ 
which satisfy the triangle inequality $f(x,\xi_1+\xi_2)\le f(x,\xi_1)+ f(x,\xi_2)$ for every 
$x\in X$ and every $\xi_1$, $\xi_2\in \Xi$.
\end{definition}

\par
\begin{remark}\label{rem:tri+lip}
As $|f(x,\xi)|\le |\xi|\,\|f\|_\hom$, each $f\in C^\hom_\triangle (X{\times}\Xi)$ is Lipschitz 
continuous with respect to $\xi$ and satisfies
$$
|f(x,\xi_1)-f(x,\xi_2)|\le |\xi_1-\xi_2|\, \|f\|_\hom
$$
 for every $x\in X$ and every $\xi_1$, $\xi_2\in \Xi$. Therefore 
$C^\hom_\triangle(X{\times}\Xi)\subset C^\hom_L(X{\times}\Xi)$.
\end{remark}

\par
\begin{lemma}\label{convhom}
The space of functions of the form $f_1-f_2$, with $f_1,\,f_2\in 
C^\hom_\triangle(X{\times}\Xi)$, is dense in $C^\hom(X{\times}\Xi)$.
\end{lemma}

\par
\begin{proof}
Thanks to the obvious density in $C^\hom(X{\times}\Xi)$ of the space of $f\in 
C^\hom(X{\times}\Xi)$ such that $f(x,\cdot)$ belongs to $C^2({\Xi\setminus\{0\}})$ 
for every $x$, it is enough to prove that every such function can be written as 
$f=f_1-f_2$, with $f_1,\,f_2\in C^\hom_\triangle(X{\times}\Xi)$. To this aim it suffices to show 
that there exists a constant $c:=c(f)$ such that $f_2(x,\xi):=c\,|\xi|-f(x,\xi)$ is convex in 
$\xi$ for every $x\in X$. A simple calculation shows that the quadratic form 
corresponding to the Hessian matrix of $f_2$ with respect to $\xi$ at a point $(x,e)$, 
with $e\in \Sigma_\Xi$, is given by
\begin{equation}\label{hessian}
D^2_\xi f_2(x,e)\,\xi{\,\cdot\,}\xi = c\,|\xi|^2-c\,(\xi{\,\cdot\,}e)^2-D^2_\xi 
f(x,e)\,\xi{\,\cdot\,}\xi\,.
\end{equation}
By the Euler relation we have $D_\xi f(x,\xi)\,\xi= f(x,\xi)$. Taking the derivative with 
respect to $\xi$ we obtain $D^2_\xi f(x,\xi)\,\xi=0$ for every $\xi$, in particular 
$D^2_\xi f(x,e)$ has an eigenvalue $0$ with eigenvector $e$. This implies that there is a 
constant $b(x,e)$ such that $D^2_\xi f(x,e) \,\xi{\,\cdot\,}\xi\le b(x,e)\, |\xi_e^\perp|^2$, 
where $\xi_e^\perp:=\xi- (\xi{\,\cdot\,}e)\,e$ is the component of $\xi$ orthogonal 
to~$e$. As $b(x,e)$ is bounded by the continuity of the second derivatives of $f$, and 
$|\xi_e^\perp|^2=|\xi|^2-(\xi{\,\cdot\,}e)^2$, by (\ref{hessian}) there exists a constant 
$c$ such that $D^2_\xi f_2(x,e)$ is positive definite for every $x\in X$ and every $e\in 
\Sigma_\Xi$, hence $f_2(x,\xi)$ is convex with respect to~$\xi$ for every $x\in X$.
\end{proof}

\par
\begin{definition}\label{M*}
The dual of the space $C^\hom(X{\times}\Xi)$ is denoted by $M_*(X{\times}\Xi)$, and 
the corresponding dual norm by $\|\cdot\|_*$; the weak$^*$ topology of 
$M_*(X{\times}\Xi)$ is defined by using this duality.
It is sometimes convenient to write the dummy variables explicitly and to use the notation 
$\langle f(x,\xi),\mu(x,\xi)\rangle$ for the duality product $\langle f,\mu\rangle$. The 
positive cone $M_*^+(X {\times}\Xi)$ is defined as the set of all 
$\mu\in M_*(X{\times}\Xi)$ such that
$$
 \langle f,\mu\rangle\ge 0  \quad \text{for every } f\in
C^\hom(X{\times}\Xi) \text{ with }
f\ge 0\,.
$$
\end{definition}

\par
\begin{remark}\label{M+}
It is easy to see that for every $\mu\in M_*^+(X {\times}\Xi)$ we have
$$
\|\mu\|_*=\langle |\xi|,\mu(x,\xi)\rangle\,.
$$
\end{remark}

\par
Strictly speaking, the elements $\mu$ of $M_*(X{\times}\Xi)$ are not measures, because 
they act only on homogeneous functions. However, the notion of image of $\mu$ under a 
map $\psi$ can be defined by duality, as in measure theory.

\par
\begin{definition}\label{def:image}
Let $\Xi$ and $\Xi'$ be two finite dimensional Hilbert spaces and let
$\psi\colon X{\times}\Xi\to X{\times}\Xi'$ be a continuous map of the form
$\psi(x,\xi)=(x,\varphi(x,\xi))$, with $\varphi\colon X{\times}\Xi\to\Xi'$ positively one-
homogeneous in $\xi$. The {\em image} $\psi(\mu)$ of $\mu\in M_*(X {\times}\Xi)$
under $\psi$ is defined as the element of $M_*(X {\times}\Xi')$ such that
$$
\langle f,\psi(\mu) \rangle =\langle f\circ\psi  ,\mu\rangle
=\langle f(x,\varphi(x,\xi)),\mu(x,\xi)\rangle
$$
for every $f\in C^\hom(X{\times}\Xi')$.
\end{definition}

\par
Similarly we can define the notion of support of $\mu\in M_*(X {\times}\Xi)$.
We say that a subset $C$ of $X{\times}\Xi$ is a $\Xi$-cone if $(x,\xi)\in C \Rightarrow 
(x,t\xi)\in C$ for every $t\ge 0$.

\par
\begin{definition}\label{def:supp}
The support $\supp\,\mu$ of $\mu\in M_*(X {\times}\Xi)$ is defined as the smallest 
closed $\Xi$-cone $C\subset X {\times}\Xi$ such that $\langle f,\mu\rangle=0$ for every 
$f\in C^\hom(X{\times}\Xi)$ vanishing on~$C$.
\end{definition}

\par
\begin{remark}\label{mutilde}
For every $\mu\in M_*(X {\times}\Xi)$ there exists a measure $\tilde\mu\in 
M_b(X{\times}\Xi)$ with compact support such that
\begin{equation}\label{integ}
\langle f,\mu\rangle=\int_{X {\times}\Xi}f\,d\tilde\mu
\end{equation}
for every $f\in C^\hom(X{\times}\Xi)$. 
A measure with this property can be constructed by considering the continuous linear map 
on $C(X{\times}\Sigma_\Xi)$ defined by $g\mapsto \langle 
|\xi|g(x,\xi/|\xi|),\mu(x,\xi)\rangle$. By the Riesz Representation Theorem there exists 
$\tilde\mu\in M_b(X{\times}\Sigma_\Xi)$ such that
$$
\langle |\xi|g(x,\xi/|\xi|),\mu(x,\xi)\rangle=\int_{X {\times}\Sigma_\Xi}g\,d\tilde\mu
$$
for every $g\in C(X{\times}\Sigma_\Xi)$.
Regarding $\tilde\mu$ as a measure on $X{\times}\Xi$ supported by 
$X{\times}\Sigma_\Xi$, we obtain (\ref{integ}). This construction suggests that the 
measure $\tilde\mu$ satisfying (\ref{integ}) is not unique; indeed, we can repeat the same 
construction with $\Sigma_\Xi$ replaced by any other concentric sphere.
\end{remark}

\par
For the applications it is convenient to extend some of the previous results to a suitable 
space of possibly discontinuous functions. 

\par
\begin{definition}\label{BXXi}
Given two finite dimensional Hilbert spaces $\Xi$ and $\Xi'$, let 
$B^\hom_\infty(X{\times}\Xi;\Xi')$ be the space of Borel functions $f\colon X 
{\times}\Xi\to \Xi'$ such that
\begin{itemize}
\item[(a)] for every $x\in X$ the function $\xi\mapsto f(x,\xi)$ is positively homogeneous 
of degree one on $\Xi$,
\item[(b)] there exists a constant $a\in \R$ such that
$|f(x,\xi)|\le a|\xi|$ for every $(x,\xi)\in X {\times}\Xi$.
\end{itemize}
The smallest constant $a$ satisfying the previous inequality is denoted by $\|f\|_{\hom}$.
When $\Xi'=\R$, the corresponding space will be denoted simply by 
$B^\hom_\infty(X{\times}\Xi)$.
\end{definition}

\par
\begin{definition}\label{fmu}
For every $f\in B^\hom_\infty(X{\times}\Xi)$ and every $\mu\in M_*(X{\times}\Xi)$ 
the duality product $\langle f,\mu\rangle$ is defined by
$$
\langle f,\mu\rangle:=\int_{X {\times}\Xi}f\,d\tilde\mu\,,
$$
where $\tilde\mu$ is any measure satisfying the conditions of Remark~\ref{mutilde}. By 
homogeneity the value of $\langle f,\mu\rangle$ does not depend on the particular 
measure $\tilde\mu$ chosen in~ (\ref{integ}).
The same definition (with values in $\R\cup\{+\infty\}$ this time) is adopted if $\mu\in 
M^+_*(X{\times}\Xi)$ and $f\colon X{\times}\Xi\to \R\cup\{+\infty\}$ is a Borel 
function such that $f(x,\xi)$ is positively homogeneous of degree one in~$\xi$ and 
$f(x,\xi)\ge -c|\xi|$ for some constant $c\ge 0$.
\end{definition}

\par
Let $\pi_X\colon X{\times}\Xi\to X$ be the projection onto $X$. We now define the 
image under $\pi_X$ of the product $h\mu$ of an element $\mu$ of $M_*(X{\times}\Xi)$ by a homogeneous function $h$.

\par
\begin{definition}\label{def:pihmu}
Let $\Xi$ and $\Xi'$ be two finite dimensional Hilbert spaces, let 
$\mu\in M_*(X{\times}\Xi)$, and let $h\in B^\hom_\infty(X{\times}\Xi;\Xi')$. The measure 
$\pi_X(h\mu)$ is the element of $M_b(X;\Xi')$ such that
\begin{equation}\label{def-hmu}
\int_X \varphi{\,\cdot\,} d\pi_X(h\mu) =\langle \varphi(x){\,\cdot\,} 
h(x,\xi),\mu(x,\xi)\rangle
\end{equation}
for every $\varphi\in C(X;\Xi')$, where the dot denotes the scalar product in $\Xi'$.
\end{definition}

\par
\begin{remark}\label{rem:piXhmu}
For every $\tilde\mu\in M_b(X{\times}\Xi)$ satisfying \eqref{integ} we can consider the 
Radon measure $h\tilde\mu\in M_b(X{\times}\Xi;\Xi')$ having density $h$ with respect 
to $\tilde\mu$. It is easy to check that the measure $\pi_X(h\mu)$ defined by 
\eqref{def-hmu} coincides with the image under $\pi_X$ of the measure $h\tilde\mu$. Note that 
the measure $h\tilde\mu$ depends on the choice of $\tilde\mu$ satisfying \eqref{integ}, 
while, by \eqref{def-hmu}, its projection $\pi_X(h\tilde\mu)$ does not. 
\end{remark}

\par
\begin{remark}\label{piXhmu}
It follows from the definition that we have the estimate
$$
\|\pi_X(h\mu)\|\le \|h\|_\hom\|\mu\|_*
$$
for every $\mu\in M_*(X {\times}\Xi)$ and every $h\in 
B^\hom_\infty(X{\times}\Xi;\Xi')$.
\end{remark}

\end{section}

\begin{section}{Generalized Young measures}

As mentioned in the introduction, the notion of generalized Young measure is used to 
describe oscillation and concentration phenomena for sequences which are bounded in 
$L^r(X;\Xi)$ only for $r=1$. To study concentration phenomena, where the sequences 
tend to infinity along given directions in the space $\Xi$, it is useful to introduce 
homogeneous coordinates. This is done by replacing the space $\Xi$ by $\Xi{\times}\R$, 
whose generic point is denoted by $(\xi,\eta)$; the set of points with $\eta=1$ is identified 
with $\Xi$, while points with $\eta=0$ are interpreted as directions at infinity.

\par
In our presentation the space of generalized Young measures will be a subset of the space
$M_*^+(X {\times}\Xi{\times} \R)$, where $\Xi{\times} \R$ plays the role of the 
Hilbert space $\Xi$ of the previous section. Before describing this set, we first consider  
generalized Young measures associated with functions.

\par
\begin{definition}\label{GYMu}
Given $u\in L^1(X;\Xi)$, the {\em generalized Young measure associated with\/} $u$ is 
defined as the element $\delta_u$ of  $M_*^+(X {\times}\Xi{\times} \R)$ such that
$$
\langle f, \delta_u \rangle = 
\int_X f(x, u(x), 1)\,d\lambda(x)
$$
for every $f\in C^\hom(X{\times}\Xi{\times}\R)$. 
\end{definition}

\par
In the spirit of \cite{Gof-Ser} and \cite{Tem} we extend this definition to measures 
$p\in M_b(X;\Xi)$.

\par
\begin{definition}\label{GYMp}
Given $p\in M_b(X;\Xi)$, the {\em generalized Young measure associated with\/} $p$ is
defined as the element $\delta_p$ of $M_*^+(X {\times}\Xi{\times} \R)$ such that for
every $f\in C^\hom(X{\times}\Xi{\times}\R)$
$$
\langle f, \delta_p \rangle = 
\int_X f(x, \textstyle \frac{dp}{d\sigma}(x),\frac{d\lambda}{d\sigma}(x))\,d\sigma(x)\,,
$$
where $\sigma$ is an arbitrary nonnegative Radon measure on $X$ with 
$\lambda<<\sigma$ and $p<<\sigma$. 
\end{definition}

\par
The homogeneity of $f$ implies that the integral does not depend on $\sigma$ and that 
the definitions coincide when $p=u\in L^1(X;\Xi)$.  The norm of $\delta_p$ is given by 
the following lemma.

\par
\begin{lemma}\label{lemma:mass}
Let $p\in M_b(X;\Xi)$. Then
$$
\|\delta_p\|_* = \int_X\sqrt{1+|p^a|^2}\, d\lambda + |p^s|(X)
\le \lambda(X)+|p|(X) \,.
$$
\end{lemma}

\par
\begin{proof}
Let us consider the Borel partition $X=X^a\cup X^s$ with 
$\lambda(X^s)=0=|p^s|(X^a)$ and let $\sigma:=\lambda +|p^s|$, so that 
$\sigma=\lambda$ on $X^a$ and $\sigma=|p^s|$ on $X^s$. By Remark~\ref{M+} we 
have
$$
\begin{array}{c}
\displaystyle
\|\delta_p\|_* = \langle \sqrt{|\xi|^2+|\eta|^2},\delta_p(x,\xi,\eta)\rangle
= \int_X \sqrt{\textstyle  | \frac{dp}{d\sigma}|^2 + |\frac{d\lambda}{d\sigma}|^2 }\, 
d\sigma =
\vspace{.2cm}
\\
\displaystyle
= \int_{X^a}\sqrt{1+|p^a|^2}\, d\lambda + |p^s|(X^s)
= \int_X\sqrt{1+|p^a|^2}\, d\lambda + |p^s|(X)\,,
\end{array}
$$
which concludes the proof.
\end{proof}

\par
We recall the definition of Young measure.

\par
\begin{definition}\label{Yr}
A {\em Young measure on $X$ with values in\/} $\Xi$ is a measure $\nu\in 
M_b^+(X{\times}\Xi)$ such that $\pi_X(\nu)=\lambda$. The space of Young measures 
on $X$ with values in $\Xi$ is denoted by $Y(X;\Xi)$. 
For every $r\ge 1$ let $Y^r(X;\Xi)$ be the space of all $\nu\in Y(X;\Xi)$ whose {\em 
$r$-moment\/}
$$
\int_{X{\times}\Xi} |\xi|^r\,d\nu(x,\xi)
$$
is finite.
\end{definition}

\par
\begin{remark}\label{Disintegration}
By the Disintegration Theorem (see, e.g., \cite[Appendix A2]{Val}) for every $\nu\in 
Y(X;\Xi)$ there exists a measurable family $(\nu^x)_{x\in X}$ of probability measures 
on $\Xi$ such that 
$$
\int_{X{\times}\Xi} g(x,\xi)\,d\nu(x,\xi)=
\int_X\Big(\int_\Xi g(x,\xi)\,d\nu^x(\xi)\Big)\, d\lambda(x)
$$
for every bounded Borel function $g\colon X{\times}\Xi\to\R$. The probability measures 
$\nu^x$ are uniquely determined for $\lambda$-a.e.\ $x\in X$.
\end{remark}

\par
\begin{definition}\label{GYMY}
Given $\nu\in Y^1(X;\Xi)$, the {\em generalized Young measure associated with\/} $\nu$ 
is defined as the element $\ol\nu$ of $M_*^+(X {\times}\Xi{\times} \R)$ such that 
$$
\langle f, \ol\nu \rangle := 
\int_{X{\times}\Xi} f(x, \xi, 1)\,d\nu(x,\xi)
$$
for every $f\in C^\hom(X{\times}\Xi{\times}\R)$. 
\end{definition}

\par
\begin{remark}\label{massYm}
It follows from Remark~\ref{M+} that
$$
\|\ol\nu\|_* = \int_{X{\times}\Xi} \sqrt{1+ |\xi|^2}\,d\nu(x,\xi)
\le \lambda(X)+ \int_{X{\times}\Xi} |\xi|\,d\nu(x,\xi) \,.
$$
\end{remark}

\par
\begin{remark}\label{properties}
If $\mu=\delta_p$ for some $p\in M_b(X;\Xi)$, the following properties hold:
\begin{eqnarray}
& \supp\,\mu\subset X {\times}\Xi{\times} [0,+\infty)\,, \label{suppmu}
\\
& \pi_X(\eta\mu) =\lambda \,.
\label{muproj}
\end{eqnarray}
We will refer to (\ref{muproj}) as the {\em projection property\/}. According to 
\eqref{def-hmu}, it is equivalent to
\begin{equation}\label{muproj2}
\langle\varphi(x)\eta, \mu(x,\xi,\eta)\rangle=\int_X\varphi \, d\lambda \quad \text{for every } 
\varphi\in
C(X)\,.
\end{equation}
Properties (\ref{suppmu}) and (\ref{muproj}) continue to hold if $\mu=\ol\nu$ for some 
$\nu\in Y^1(X;\Xi)$. 
\end{remark}

\par
This motivates the following definition.

\par
\begin{definition}\label{GY}
The space $GY(X;\Xi)$ of {\it generalized Young measures\/} on $X$ with values in 
$\Xi$ is defined as the set of all $\mu\in M_*^+(X {\times}\Xi{\times} \R)$ satisfying 
(\ref{suppmu}) and (\ref{muproj}). On $GY(X;\Xi)$ we consider the norm and the 
weak$^*$ topology induced by $M_*(X {\times}\Xi{\times} \R)$.
\end{definition}

\par
\begin{remark}\label{boundedborel}
By approximation we can prove that (\ref{muproj2}) holds for every $\mu\in GY(X;\Xi)$ 
and for every bounded Borel function $\varphi\colon X\to\R$.
\end{remark}

\par
The sequential compactness of every bounded subset of $GY(X;\Xi)$ is given by the 
following theorem.

\par
\begin{theorem}\label{thm:cpt}
Every bounded sequence in $GY(X;\Xi)$ has a subsequence which converges  
weakly$^*$ to an element of $GY(X;\Xi)$.
\end{theorem}

\par
\begin{proof}
Since $GY(X;\Xi)$ is closed in the weak$^*$ topology of $M_*(X {\times}\Xi{\times} 
\R)$, 
the result follows from the Banach-Alaoglu Theorem.
\end{proof}

\par
\begin{remark}\label{convnorms}
If $\mu_k$ is a sequence in $GY(X;\Xi)$ which converges weakly$^*$ to $\mu\in 
GY(X;\Xi)$, then $\|\mu_k\|_*\to \|\mu\|_*$ by Remark~\ref{M+}.
\end{remark}

\par
\begin{remark}\label{lscw*}
If $f\colon X{\times}\Xi{\times}\R\to \R\cup\{+\infty\}$ is lower semicontinuous, 
$(\xi,\eta)\mapsto f(x,\xi,\eta)$ is positively homogeneous of degree one, and 
$f(x.\xi,\eta)\ge-c\sqrt{|\xi|^2+|\eta|^2}$ for some constant $c\ge 0$, then
$$
\langle f,\mu\rangle \le \liminf_{k\to\infty}\, \langle f,\mu_k\rangle
$$
for every sequence $\mu_k$ in $GY(X;\Xi)$ which converges weakly$^*$ to $\mu\in 
GY(X;\Xi)$, where $\langle f,\cdot\rangle$ is defined by~\eqref{integ}. Indeed, any such 
$f$ is the supremum of a family of functions in $C^\hom(X{\times}\Xi{\times}\R)$.
\end{remark}

\par
In the case of a generalized Young measure $\mu\in GY(X;\Xi)$  the duality product 
$\langle f,\mu\rangle$ can be defined for every $f$ in the space 
$B^\hom_{\infty,1}(X{\times}\Xi{\times}\R)$ introduced by the following definition, 
which is slightly larger than the space $B^\hom_\infty(X{\times}\Xi{\times}\R)$ 
considered in the previous section.

\par
\begin{definition}\label{B1XXi}
Given two finite dimensional Hilbert spaces $\Xi$ and $\Xi'$, we consider the space 
$B^\hom_{\infty,1}(X{\times}\Xi{\times}\R;\Xi')$ of all Borel functions $f\colon X 
{\times}\Xi{\times}\R\to \Xi'$ such that
\begin{itemize}
\item[(a)] for every $x\in X$ the function $(\xi,\eta)\mapsto f(x,\xi,\eta)$ is positively 
homogeneous of degree one on $\Xi{\times}\R$,
\item[(b)] there exist a constant $a\in \R$ and a function $b\in L^1(X)$ such that
\begin{equation}\label{ling}
|f(x,\xi,\eta)|\le a|\xi|+b(x)|\eta|
\end{equation}
for every $(x,\xi,\eta)\in X {\times}\Xi{\times}\R$.
\end{itemize}
When $\Xi'=\R$, the corresponding space will be denoted simply by 
$B^\hom_{\infty,1}(X{\times}\Xi{\times}\R)$.
\end{definition}

\par
\begin{lemma}\label{fL1}
Let $\mu\in GY(X;\Xi)$ and let $\tilde\mu\in M_b(X{\times}\Xi{\times}\R)$ be a 
measure with compact support satisfying (\ref{integ}). Then every $f\in 
B^\hom_{\infty,1}(X{\times}\Xi{\times}\R)$ is $\tilde\mu$-integrable. 
\end{lemma}

\par
\begin{proof}
Let us fix $f\in B^\hom_{\infty,1}(X{\times}\Xi{\times}\R)$.
For every $k$ let $f_k$ be defined by $f_k(x,\xi,\eta):=f(x,\xi,\eta)$ if $|f(x,\xi,\eta)|\le 
a|\xi|+ k|\eta|$, and by  $f_k(x,\xi,\eta):=0$ otherwise. Then we have 
\begin{eqnarray*}
& \displaystyle
\int_{X{\times}\Xi{\times}\R} |f_k|\,d\tilde\mu \le \int_{X{\times}\Xi{\times}\R} [a|\xi|+ 
(b(x)\land k)|\eta|]\,d\tilde\mu = \langle a|\xi|+(b(x)\land k) |\eta|, \mu(x,\xi,\eta)\rangle \le
\\
& \displaystyle
\le \langle a|\xi|, \mu(x,\xi,\eta)\rangle +\int_X b\, d\lambda \,.
\end{eqnarray*}
It follows from Fatou's Lemma that $f$ is $\tilde\mu$-integrable.
\end{proof}

\par
\begin{definition} Given $f\in B^\hom_{\infty,1}(X{\times}\Xi{\times}\R)$ and $\mu\in 
GY(X;\Xi)$,
the duality product $\langle f,\mu\rangle$ is defined by 
\begin{equation}\label{integ3}
\langle f,\mu\rangle:=\int_{X{\times}\Xi{\times}\R} f\,d\tilde\mu\,,
\end{equation}
where $\tilde\mu\in M_b(X{\times}\Xi{\times}\R)$ is any measure with compact support 
satisfying (\ref{integ}), with $\Xi$ replaced by~$\Xi{\times}\R$.
\end{definition}

\par
\begin{remark}\label{rem:fmu10}
The integral in \eqref{integ3} is well defined by Lemma~\ref{fL1}.
It is easy to see that the value of this integral does not depend on the choice of 
$\tilde\mu$ 
satisfying (\ref{integ}), and that
$$
|\langle f,\mu\rangle | \le a\,\|\mu\|_* +\| b\|_1\,,
$$
where $a$ and $b$ satisfy \eqref{ling} and $\| b\|_1$ denotes the $L^1$ norm of $b$.
\end{remark}

\par
We now consider the image of a generalized Young measure.

\par
\begin{definition}\label{def:image3}
Let $\Xi$ and $\Xi'$ be two finite dimensional Hilbert spaces and let $\psi\colon 
\!X{\times}\Xi{\times}\R\allowbreak\to X{\times}\Xi'{\times}\R $ be a map of the form 
$\psi(x,\xi,\eta)\!=\!(x,\varphi(x,\xi,\eta),\eta)$, with $\varphi\!\in 
\!B^\hom_{\infty,1}(X{\times}\Xi{\times}\R;\Xi')$. The {\em image} $\psi(\mu)$ of 
$\mu\in GY(X;\Xi)$ under $\psi$ is defined as the element of $GY(X;\Xi')$ such that 
\begin{equation}\label{image3}
\langle f,\psi(\mu) \rangle =\langle f\circ\psi  ,\mu\rangle
=\langle f(x,\varphi(x,\xi,\eta),\eta),\mu(x,\xi,\eta)\rangle
\end{equation}
for every $f\in B^\hom_{\infty,1}(X{\times}\Xi'{\times}\R)$.
\end{definition}

\par
\begin{remark}\label{rem:image3}
Under the assumptions of the definition the function $f\circ\psi$ belongs to 
$B^\hom_{\infty,1}(X{\times}\Xi{\times}\R)$, so that the duality product $\langle 
f\circ\psi  ,\mu\rangle$ is well defined. Moreover, by the particular form of the map $\psi$, 
the element of $M^+_*(X{\times}\Xi'{\times}\R)$ defined by \eqref{image3} satisfies 
\eqref{suppmu} and \eqref{muproj}, therefore it belongs to $GY(X;\Xi')$.
\end{remark}

\end{section}

\begin{section}{Comparison with other presentations of the theory}

In this section we show that every $\mu\in GY(X;\Xi)$ can be represented by a unique 
Young measure-varifold pair $(\mu^Y\!,\mu^\infty)$, where $\mu^Y\!\in Y^1(X;\Xi)$ and 
$\mu^\infty\in M_b^+(X{\times}\Sigma_\Xi)$. To introduce this representation, we recall 
that $Y^1(X;\Xi)$ can be identified with a suitable subset of $GY(X;\Xi)$ 
(Definition~\ref{GYMY}). The following definition identifies the measures in 
$M_b^+(X{\times}\Sigma_\Xi)$ with particular elements of 
$M_*(X{\times}\Xi{\times}\R)$.

\par
\begin{definition}\label{hatnu}
For every $\nu\in M_b^+(X{\times}\Sigma_\Xi)$ let $\hat\nu$ be the element of 
$M_*(X{\times}\Xi{\times}\R)$ defined by
$$
\langle f,\hat\nu\rangle =\int_{X{\times}\Sigma_\Xi}f(x,\xi,0)\,d\nu(x,\xi)
$$
for every $f\in C^\hom(X{\times}\Xi{\times}\R)$. 
\end{definition}

\par
\begin{remark}\label{rmk:145}
It follows from the definition that $\pi_X(\eta\hat\nu)=0$ for every $\nu\in 
M_b^+(X{\times}\Sigma_\Xi)$. Since $\hat\nu$ does not satisfy the projection property 
\eqref{muproj}, it does not belong to $GY(X;\Xi)$.
\end{remark}

\par
The main result of this section is the following theorem.

\par
\begin{theorem}\label{thm:decomp}
Let $\mu\in GY(X;\Xi)$. Then there exists a unique pair $(\mu^Y\!, \mu^\infty)$, with 
$\mu^Y\!\in Y^1(X;\Xi)$ and $\mu^\infty\in  M_b^+(X{\times}\Sigma_\Xi)$, such that 
$$
\mu=\ol\mu^Y+\hat\mu^\infty\,,
$$
which is equivalent to
\begin{equation}\label{mu1muinf}
\langle f,\mu\rangle = \int_{X{\times}\Xi} f(x,\xi,1)\,d\mu^Y\!(x,\xi) +
\int_{X{\times}\Sigma_\Xi}f(x,\xi,0)\,d\mu^\infty(x,\xi)
\end{equation}
for every $f\in B^\hom_{\infty,1}(X{\times}\Xi{\times}\R)$. 
\end{theorem}

\par
\begin{remark}\label{rmk:1}
The converse of Theorem~\ref{thm:decomp} is also true: if $\mu^Y\!\in Y^1(X;\Xi)$ and 
$\mu^\infty\in M_b(X{\times}\Sigma_\Xi)$, then formula (\ref{mu1muinf}) defines an 
element of $GY(X;\Xi)$.
\end{remark}

\par
\begin{remark}\label{disint}
Let $\lambda^\infty:=\pi_X(\mu^\infty)$. Since $\lambda=\pi_X(\mu^Y)$, by the 
Disintegration Theorem (see, e.g., \cite[Appendix A2]{Val}) there exist a measurable 
family $(\mu^{x,Y})_{x\in X}$ of probability measures on $\Xi$ and a measurable 
family $(\mu^{x,\infty})_{x\in X}$ of probability measures on $\Sigma_\Xi$ such that 
$$
\begin{array}{c}
\displaystyle \langle f,\mu\rangle = \int_{X{\times}\Xi} f(x,\xi,1)\,d\mu^Y(x,\xi) +
\int_{X{\times}\Sigma_\Xi}f(x,\xi,0)\,d\mu^\infty(x,\xi) =
\vspace{.2cm}
\\
\displaystyle = \int_X\Big(\int_\Xi f(x,\xi,1)\,d\mu^{x,Y}(\xi)\Big)\, d\lambda(x)+
\int_X\Big(\int_{\Sigma_\Xi} f(x,\xi,0)\,d\mu^{x,\infty}(\xi)\Big)\, d\lambda^\infty(x)
\end{array}
$$
for every $f\in B^\hom_{\infty,1}(X{\times}\Xi{\times}\R)$. 
\end{remark}

\par
\begin{remark}\label{rmk:norm}
Thanks to Remark~\ref{M+}, if we apply (\ref{mu1muinf}) to 
$f(x,\xi,\eta)=\sqrt{|\xi|^2+|\eta|^2}$, we obtain
$$
\begin{array}{c}
\displaystyle
\|\mu\|_*= \int_{X{\times}\Xi} \sqrt{1+|\xi|^2} \,d\mu^Y(x,\xi)  + 
\mu^\infty(X{\times}\Sigma_\Xi)
\le 
\\
\displaystyle \le \lambda(X) + \int_{X{\times}\Xi}|\xi|\,d\mu^Y(x,\xi)  + 
\mu^\infty(X{\times}\Sigma_\Xi)\,.
\end{array}
$$
\end{remark}

\par
\begin{proof}[Proof of Theorem~\ref{thm:decomp}]
For every Borel function $g\colon X{\times}\Xi\to\R$ with
\begin{equation}\label{gLin}
\kappa_g:=\sup_{(x,\xi)\in X{\times}\Xi}\frac{|g(x,\xi)|}{\sqrt{1+|\xi|^2}}<+\infty\,,
\end{equation}
we consider the Borel function $\varphi_g\colon X{\times}\Xi{\times}\R\to\R$ defined 
by
$$
\varphi_g(x,\xi,\eta): =
\begin{cases}
\eta\, g(x,\xi/\eta)&
\text{if } \eta>0\,,
\\
0&
\text{if } \eta\le0\,.
\end{cases}
$$
Since $\varphi_g\in B^\hom_\infty(X{\times}\Xi{\times}\R)$, we can consider the duality 
product $\langle \varphi_g,\mu\rangle$. The function $g\mapsto \langle 
\varphi_g,\mu\rangle$ is linear, bounded, and positive on $C_0(X{\times}\Xi)$. By the 
Riesz Representation Theorem there exists $\mu^Y\!\in M_b^+(X{\times}\Xi)$ such that
\begin{equation}\label{defmu1}
\langle \varphi_g,\mu\rangle = \int_{X{\times}\Xi} g(x,\xi)\,d\mu^Y\!(x,\xi)
\end{equation}
for every $g\in C_0(X{\times}\Xi)$. As $\|g\|_\hom\le \kappa_g$, we have
$$
|\int_{X{\times}\Xi} g(x,\xi)\,d\mu^Y\!(x,\xi)| \le  \kappa_g\|\mu\|_*
$$
for every $g\in C_0(X{\times}\Xi)$. By approximation we can prove that 
\begin{equation}\label{mom122}
\int_{X{\times}\Xi} |\xi|\,d\mu^Y\!(x,\xi)\le \|\mu\|_*
\end{equation}
and that \eqref{defmu1} holds for every Borel function $g\colon X{\times}\Xi\to\R$ 
satisfying (\ref{gLin}).
By (\ref{muproj}) we have $\pi_{X}(\mu^Y)=\lambda$,
which, together with \eqref{mom122}, gives $\mu^Y\!\in Y^1(X;\Xi)$. 

\par
For every bounded Borel function $h\colon X{\times}\Sigma_\Xi\to\R$ we consider the 
Borel function 
$\psi_h\colon \allowbreak X{\times}\Xi{\times}\R\to\R$ defined by
$$
\psi_h(x,\xi,\eta): =
\begin{cases}
|\xi|\, h(x,\xi/|\xi|)&
\text{if } \eta=0\,,\xi\neq 0\,,
\\
0&
\text{otherwise.}
\end{cases}
$$
Since $\psi_h\in B^\hom_\infty(X{\times}\Xi{\times}\R)$, we can consider the duality 
product $\langle \psi_h,\mu\rangle$. The function $h\mapsto \langle \psi_h,\mu\rangle$ is  
linear, bounded, and positive  on $C(X{\times}\Sigma_\Xi)$. By the Riesz Representation 
Theorem there exists $\mu^\infty\in M_b^+(X{\times}\Sigma_\Xi)$ such that
\begin{equation}\label{defmuinfty}
\langle \psi_h,\mu\rangle = \int_{X{\times}\Sigma_\Xi} h\,d\mu^\infty
\end{equation}
for every $h\in C(X{\times}\Sigma_\Xi)$. By approximation we can prove that the 
previous equality holds for every bounded Borel function $h\colon 
X{\times}\Sigma_\Xi\to\R$.

\par
Given any $f\in B^\hom_\infty(X{\times}\Xi{\times}\R)$, we consider the functions 
$g\colon\! X{\times}\Xi\to\R$ and $h\colon\! X{\times}\Sigma_\Xi\to\R$ defined by
$$
g(x,\xi):=f(x,\xi,1)\,, \qquad h(x,\xi):=f(x,\xi,0)\,.
$$
By homogeneity we have $f=\varphi_g+\psi_h$ on $X{\times}\Xi{\times}[0,+\infty)$. 
Then (\ref{mu1muinf}) follows from \eqref{suppmu}, (\ref{defmu1}), and 
(\ref{defmuinfty}). The result can be extended to $f\in 
B^{\hom,1}_\infty(X{\times}\Xi{\times}\R)$ by approximation.

\par
The uniqueness of the pair $(\mu^Y\!,\mu^\infty)$ can be deduced from the fact that, 
if \eqref{mu1muinf} is satisfied, then \eqref{defmu1} holds for every 
$g\in C_0(X{\times}\Xi)$, while \eqref{defmuinfty} holds for every 
$h\in C(X{\times}\Sigma_\Xi)$.
\end{proof}

\par
\begin{remark}\label{momenta}
It is easy to prove, by approximation, that \eqref{mu1muinf} continues to hold when 
$f\colon X{\times}\Xi{\times}\R\to \R\cup\{+\infty\}$ is a Borel function such that 
$f(x,\xi,\eta)$ is positively homogeneous of degree one in~$(\xi,\eta)$ and satisfies the 
inequality $f(x,\xi,\eta)\ge -c\sqrt{|\xi|^2+|\eta|^2}$ for some constant $c\ge 0$. This 
allows to characterize the generalized Young measures associated with Young measures 
with finite $r$-moment, $r>1$, using the homogeneous functions
$\{P_r\}\colon\Xi{\times}\R\to[0,+\infty]$ defined by
\begin{equation}\label{mom1}
\{P_r\}(\xi,\eta):= \begin{cases}
|\xi|^r/\eta^{r-1} & \text{if } \eta>0\,,
\\
+\infty & \text{if } \eta\le 0\,.
\end{cases}
\end{equation}
Indeed, for $r>1$ formula \eqref{mu1muinf} implies that $\mu=\ol\mu^Y$ with 
$\mu^Y\in Y^r(X;\Xi)$ if and only if $\langle \{P_r\},\mu \rangle<+\infty$.
\end{remark}

\par
Let $\psi_0^\Xi\colon X{\times}\Xi{\times}\R\to X{\times}\Xi{\times}\R$ be the Borel 
map defined by
$$
\psi_0^\Xi (x,\xi,\eta)=
\begin{cases}
(x,\xi,\eta) & \text{if } \eta\neq 0\,,
\\
(x,0,0) & \text{if } \eta=0\,.
\end{cases}
$$
Note that $\psi_0^\Xi$ satisfies the conditions of Definition~\ref{def:image3}.

\par
For every $f\in B^\hom(X{\times}\Xi{\times}\R)$ we have
$$
(f\circ\psi_0^\Xi)(x,\xi,\eta)=
\begin{cases}
f(x,\xi,\eta) & \text{if } \eta\neq 0\,,
\\
0 & \text{if } \eta=0\,.
\end{cases}
$$
From (\ref{mu1muinf}) it follows that for every $\mu\in GY(X;\Xi)$
$$
\langle f,\ol\mu^Y\rangle=\langle f\circ\psi_0^\Xi, \ol\mu^Y\rangle=\langle 
f\circ\psi_0^\Xi,\mu\rangle\,,
$$
hence $\ol\mu^Y=\psi_0^\Xi(\mu)$.

\par
\begin{lemma}\label{imagefin}
Let $\Xi$ and $\Xi'$ be two finite dimensional Hilbert spaces, let $\mu\in GY(X;\Xi)$, let 
$\psi\colon X{\times}\Xi{\times}\R\to X{\times}\Xi'{\times}\R$ be a map as in 
Definition~\ref{def:image3}, and let $\nu:=\psi(\mu)$. Then 
$$
\ol\nu^Y=\psi(\ol\mu^Y)\,, \qquad \hat\nu^\infty=\psi(\hat\mu^\infty)\,.
$$
\end{lemma}

\par
\begin{proof}
The former equality follows from the fact that 
$\psi\circ\psi_0^\Xi=\psi_0^{\Xi'}\circ\psi$. The latter follows now from 
Theorem~\ref{thm:decomp} by the linearity of the map $\mu\mapsto \psi(\mu)$.
\end{proof}

\par
Combining the compactness property (Theorem~\ref{thm:cpt}) and the representation 
formula (Theorem~\ref{thm:decomp}) we recover the following result, originally proved
in~\cite{Ali-Bou} (see Remark~\ref{disint}).

\par
\begin{theorem}\label{th432}
Let $u_k$ be a bounded sequence in $L^1(X;\Xi)$. Then there exist a subsequence,
still denoted $u_k$, a Young measure $\mu^Y\!\in Y^1(X;\Xi)$, and a measure 
$\mu^\infty\in M_b(X{\times}\Sigma_\Xi)$, such that 
\begin{equation}
\int_X g(x,u_k(x))\,d\lambda(x) \longrightarrow
\int_{X{\times}\Xi} g(x,\xi)\, d\mu^Y\!(x,\xi) +
\int_{X{\times}\Sigma_\Xi} g^\infty (x,\xi)\, d\mu^\infty(x,\xi)
\end{equation}
for every continuous function $g\colon X{\times}\Xi\to\R$ such that for every
$(x_0,\xi_0)\in X{\times}\Xi$ the limit
$$
g^\infty(x_0,\xi_0):=
\lim_{\genfrac{}{}{0pt}2{\scriptstyle x\to x_0,\ \xi\to\xi_0}{\scriptstyle \eta\to 0^+}}
\eta\, g(x,\xi/\eta)
$$
exists and is finite.
\end{theorem}

\par
\begin{proof}
Let us consider the sequence $\delta_{u_k}$ in $GY(X;\Xi)$ introduced in
Definition~\ref{GYMu}. By Lemma~\ref{lemma:mass} we have
$\|\delta_{u_k}\|_*\le \sup_j \| u_j\|_1 +\lambda(X)<+\infty$. 
By Theorem~\ref{thm:cpt} there exists a subsequence, still denoted $u_k$, such that
$\delta_{u_k}$ converge weakly$^*$ to an element $\mu$ of $GY(X;\Xi)$.
Let $g$ be as in the statement of the theorem and let
$f\colon X{\times}\Xi{\times}\R\to\R$ be defined by
$$
f(x,\xi,\eta):=
\begin{cases}
\eta\, g(x,\xi/\eta)&
\text{if } \eta>0\,,
\\
g^\infty(x,\xi)&
\text{if } \eta\le0\,.
\end{cases}
$$
It is easy to check that $f$ is continuous in $(x,\xi,\eta)$ and homogeneous of degree one 
in $(\xi,\eta)$. Therefore, the weak$^*$ convergence of $\delta_{u_k}$ to $\mu$ implies 
that
\begin{equation}\label{alib1}
\int_X g(x,u_k(x))\,d\lambda(x)=
\int_X f(x,u_k(x),1)\,d\lambda(x) \longrightarrow
\langle f, \mu\rangle\,.
\end{equation}
By Theorem~\ref{thm:decomp}, taking into account the definition of $f$, we obtain that 
there exists a pair $(\mu^Y\!,\mu^\infty)$, with $\mu^Y\!\in Y^1(X;\Xi)$ and 
$\mu^\infty\in  M_b^+(X{\times}\Sigma_\Xi)$, such that  
\begin{equation}\label{alib2}
\langle f, \mu\rangle = \int_{X{\times}\Xi} g(x,\xi)\, d\mu^Y\!(x,\xi) +
\int_{X{\times}\Sigma_\Xi} g^\infty (x,\xi)\, d\mu^\infty(x,\xi)\,.
\end{equation}
The conclusion follows from \eqref{alib1} and \eqref{alib2}.
\end{proof}

\end{section}

\begin{section}{A density result}

In this section we prove that, if $\lambda$ is nonatomic, then the generalized Young 
measures of the form $\delta_u$ associated with functions $u\in L^1(X;\Xi)$ are dense in 
$GY(X;\Xi)$. The main result is the following approximation theorem.

\par
\begin{theorem}\label{thm:density}
Assume that $\lambda$ is nonatomic and 
let $\mu\in GY(X;\Xi)$. Then there exists a sequence $u_n$ in $L^1(X;\Xi)$ such that 
$\delta_{u_n}\wto \mu$ weakly$^*$ in $GY(X;\Xi)$.
\end{theorem}

\par
\begin{proof}
We consider the decomposition
$$
\mu=\ol\mu^Y + \hat\mu^\infty
$$
of Theorem~\ref{thm:decomp} and we fix a sequence $\sigma_n$ converging to~$0$ 
with $0<\sigma_n<\min\{1,\lambda(X)\}$.  For every $n$ we consider two countable 
partitions $\Xi=\bigcup_{j} B^{n,Y}_j$ and $\Sigma_\Xi=\bigcup_{j} B^{n,\infty}_j$, 
where the sets $B^{n,Y}_j$ and $ B^{n,\infty}_j$ satisfy
\begin{equation}\label{diam1}
{\rm diam\,} B^{n,Y}_j\le \sigma_n\qquad\hbox{and}\qquad
{\rm diam\,} B^{n,\infty}_j\le \sigma_n \,.
\end{equation}

\par
Let $\lambda^Y\!:=\pi_X(\sqrt{1+|\xi|^2}\mu^Y)$ and 
$\lambda^\infty:=\pi_X(\mu^\infty)$; i.e., 
$$
\lambda^Y\!(A)=\int_{A{\times}\Xi}\sqrt{1+|\xi|^2}\, d\mu^Y\!(x,\xi) \qquad 
\hbox{and} \qquad
\lambda^\infty(A)=\mu^\infty(A{\times}\Sigma_\Xi)
$$
for every Borel set $A\subset X$. As $\lambda=\pi_X(\mu^Y\!)$, the measure 
$\lambda^Y$ is absolutely continuous with respect to $\lambda$.

\par
Let us consider the partition $X=X^a\cup X^s$, with 
$\lambda(X^s)=0=\lambda^{\infty,s}(X^a)$, where $\lambda^{\infty,s}$ is the singular 
part of $\lambda^\infty$ with respect to $\lambda$. 
To prove the theorem, for every $n$ we will consider a new partition $X=X^{n,a}\cup 
X^{n,s}$, where $X^{n,a}$ and $X^{n,s}$  are suitable approximations of $X^a$ and 
$X^s$ such that $\lambda(X^{n,a})>0$ and $\lambda(X^{n,s})>0$. 
We will construct the approximating sequence $u_n$ by defining it separately on 
$X^{n,a}$ and $X^{n,s}$. 

\par
\medskip
\noindent
{\it Step 1. Definition of $u_n$ on $X^{n,s}$.\/}
We begin by constructing $X^{n,s}$. For every $n$ we can find a countable Borel 
partition $X^s=\bigcup_i X_i^{n,s}\cup N^{n,s}$, where each $X_i^{n,s}$ is closed, 
\begin{equation} \label{diam-2}
{\rm diam\,}X^{n,s}_i \le  \sigma_n/2 \,,\qquad\hbox{and}\qquad
\lambda^\infty(N^{n,s})=0\,. 
\end{equation}
In the following, given a subset $E\subset X$ and a radius $r>0$, the $r$-neighbourhood 
of $E$ will be denoted by 
$$
(E)_r:=\{ x\in X: \ d(x,E)<r\}\,.
$$
Since $\lambda(X_i^{n,s})= \lambda^Y\!(X_i^{n,s})=0$ and 
$\lambda((X_i^{n,s})_r)>0$ for every $r>0$, we can construct inductively a decreasing 
sequence $r^{n}_i$ such that $0<r^n_i\le\sigma_n/2$,  
\begin{eqnarray}
&
\lambda((X_i^{n,s})_{r^n_i})\le \sigma_n \lambda^\infty(X_i^{n,s})\,, \label{diam-1}\\
\smallskip
&
\lambda^Y\!((X_i^{n,s})_{r^n_i})\le 2^{-i}\sigma_n\,, \label{lambdaY}
\\
\smallskip
&
\lambda^\infty((X_i^{n,s})_{r^n_i}\setminus X^s)\le 2^{-i}\sigma_n\,,
\label{701}
\\
\smallskip
&
\lambda((X_{i+1}^{n,s})_{r^n_{i+1}})\le \frac{1}{3}\lambda((X_i^{n,s})_{r^n_i})\,.
\label{750}
\end{eqnarray}
We define
$$
\textstyle
A^{n,s}_i:= (X_i^{n,s})_{r^n_i}\setminus \bigcup_{j>i} (X_j^{n,s})_{r^n_j}\,.
$$
By \eqref{750} we have
$$
\textstyle
\lambda(A^{n,s}_i)\ge \frac{1}{2} \lambda((X_i^{n,s})_{r^n_i})>0\,,
$$
while \eqref{diam-2}, together with the inequality $0<r^n_i\le\sigma_n/2$, yields
\begin{equation}\label{diam0}
{\rm diam\,}A^{n,s}_i \le \sigma_n\,.
\end{equation}
By \eqref{diam-1} we have
$$
0<\lambda(A^{n,s}_i)\le \sigma_n\lambda^{\infty}(X_i^{n,s})=
\sigma_n\mu^{\infty}(X_i^{n,s}{\times}\Sigma_\Xi)
=\sigma_n\textstyle\sum_{j}\mu^{\infty}(X_i^{n,s}{\times}B^{n,\infty}_j)\,.
$$
Since $\lambda$ is nonatomic 
we can find a countable Borel partition $A^{n,s}_i=\bigcup_j A^{n,s}_{ij}$ such that 
\begin{equation}\label{choice0}
0<\lambda(A^{n,s}_{ij})\le \sigma_n \mu^{\infty}(X_i^{n,s}{\times}B^{n,\infty}_j)\,.
\end{equation}
For every $n$ we define 
\begin{equation}\label{704}
\textstyle X^{n,s}:= \bigcup_{ij} A^{n,s}_{ij} \cup X^s \,.
\end{equation}
Note that by \eqref{lambdaY} and \eqref{701} we have
\begin{equation}\label{702}
\lambda^Y\!(X^{n,s})\le \sigma_n \qquad \text{and} \qquad 
\lambda^\infty(X^{n,s}\setminus X^s)\le \sigma_n \,.
\end{equation} 
We define 
\begin{equation}\label{pndef4}
u_n(x):=c^{n}_{ij}\xi^{n,\infty}_j\quad\text{for }x\in A^{n,s}_{ij}\,,
\end{equation}
where $\xi^{n,\infty}_j$ are arbitrary points of $B^{n,\infty}_j$ and
\begin{equation}\label{defcnk}
c^{n}_{ij}:=\frac{\mu^\infty( 
X_i^{n,s}{\times}B^{n,\infty}_j)}{\lambda({A^{n,s}_{ij}})}\,.
\end{equation}
By \eqref{choice0} we have that 
\begin{equation}\label{choice0.5}
c^n_{ij}\ge 1/\sigma_n\,.
\end{equation}
By \eqref{defcnk} and by \eqref{choice0.5} we have 
\begin{equation}\label{check3}
\int_{X^{n,s}}\sqrt{1+|u_n|^2}\,d\lambda \le \lambda^{\infty}(X^s) 
\sqrt{1+\sigma_n^2}\,.
\end{equation}

\par
\medskip
\noindent
{\it Step 2. Definition of $u_n$ on $X^{n,a}$.\/}
We set
$$
X^{n,a}:=X\setminus X^{n,s}\,.
$$
In order to define $u_n$ on $X^{n,a}$  we consider a countable Borel partition 
$X^{n,a}=\bigcup_{i} A^{n,a}_i$, with $A^{n,a}_i$ satisfying
\begin{equation}\label{diam2}
0<{\rm diam\,}A^{n,a}_i\le \sigma_n\,.
\end{equation}
As  $X^{n,a}\subset X^a$ by \eqref{704},  $\lambda^\infty$ is absolutely continuous 
with respect to $\lambda$ on $X^{n,a}$. Since $\lambda$ is nonatomic, for every $i$ we 
may choose $0<\e^n_i\le\sigma_n$ and two disjoint Borel sets $A^{n,Y}_i\!$ and 
$A^{n,\infty}_i$ in such a way that $A^{n,a}_i=A^{n,Y}_i\!\cup A^{n,\infty}_i$ and 
\begin{equation}\label{smallportion}
\lambda(A^{n,\infty}_i)=\e^n_i\lambda^\infty(A^{n,a}_i)\le \sigma_n 
\lambda(A^{n,a}_i)\,.
\end{equation}
Since $\lambda$ is nonatomic and
$$
\e^n_i \lambda^\infty(A^{n,a}_i)=\e^n_i\mu^\infty(A^{n,a}_i{\times}\Sigma_\Xi) 
=\textstyle\sum_{j}\e^n_i\mu^\infty(A^{n,a}_i{\times}B^{n,\infty}_j)\,,
$$
we can also find a countable Borel  partition $A^{n,\infty}_{i}=\bigcup_j 
A^{n,\infty}_{ij}$  such that
\begin{equation}\label{choice1}
\lambda(A^{n,\infty}_{ij})=\e^n_i\mu^\infty(A^{n,a}_i{\times}B^{n,\infty}_j)\,.
\end{equation}
Note also that by \eqref{smallportion} we have 
$$
\lambda(A^{n,Y}_i\!)=\lambda(A^{n,a}_i)-\lambda(A^{n,\infty}_i)\ge (1-
\sigma_n)\,\lambda(A^{n,a}_i)
$$
and so there exists $0<\delta^n_i\le\sigma_n$ such that $\lambda(A^{n,Y}_i\!)=
(1-\delta^n_i)\lambda(A^{n,a}_i)$. As $\lambda=\pi_X(\mu^Y\!)$, arguing as before we 
may  find a countable Borel partition $A^{n,Y}_i\!=\bigcup_jA^{n,Y}_{ij}\!$  such that
\begin{equation}\label{choice2}
\lambda(A^{n,Y}_{ij}\!)=(1-\delta^n_i)\mu^Y\!(A^{n,a}_i{\times}B^{n,Y}_j\!)\,.
\end{equation} 

\par
We are ready to define $u_n$ on $X^{n,a}$ by setting 
\begin{equation}\label{pndef1}
u_n(x):=\xi^{n,Y}_j\quad\text{for } x\in A^{n,Y}_{ij}
\qquad\hbox{and}\qquad
u_n(x):=\frac{1}{\e^n_i}\xi^{n,\infty}_j\quad\text{for } x\in A^{n,\infty}_{ij}\,,
\end{equation}
where $\xi^{n,Y}_j\!$ are arbitrary points in $B^{n,Y}_j\!$ and $\xi^{n,\infty}_j$ are 
the points of $B^{n,\infty}_j$ chosen in~\eqref{pndef4}.
Using \eqref{choice1} and \eqref{choice2} it is easy to check that 
\begin{equation}\label{check1}
\int_{X^{n,a}}\sqrt{1+|u_n|^2}\,d\lambda\le  \sigma_n\lambda(X)+
\lambda^Y(X) +
\lambda^{\infty}(X^{n,a})\sqrt{1+\sigma_n^2}\,.
\end{equation}
By \eqref{check3} and \eqref{check1} we have
\begin{equation}\label{check2}
\int_X \sqrt{1+|u_n|^2}\,d\lambda\le  \sigma_n\lambda(X)+
\lambda^Y(X) +
\lambda^{\infty}(X)\sqrt{1+\sigma_n^2}\,,
\end{equation}
which implies that $u_n$ is bounded in $L^1(X;\Xi)$. It follows from
Lemma~\ref{lemma:mass} that $\|\delta_{u_n}\|_*$ is uniformly bounded.

\par
\medskip
\noindent
{\it Step 3. Proof of the convergence.\/}
Thanks to Lemma~\ref{L-density}, to prove the weak$^*$ convergence of 
$\delta_{u_n}$ it is enough to show that
\begin{equation}\label{499}
\langle f, \delta_{u_n}\rangle\to \langle f, \mu\rangle,
\end{equation}
for every  function $f\in C^{\hom}_L(X{\times}\Xi{\times}\R)$. Let us fix $f\in 
C^{\hom}_L(X{\times}\Xi{\times}\R)$. By Remark~\ref{rem:ChomL} there exist a 
constant $a\in\R$ and a continuous functon $\omega\colon[0+\infty)\to [0+\infty)$, with 
$\omega(0)=0$, such that
\begin{equation}\label{500}
\begin{array}{c}
|f(x_1,\xi_1,\eta_1)-f(x_2,\xi_2,\eta_2)|\leq
a \sqrt{|\xi_1-\xi_2|^2+|\eta_1-\eta_2|^2} +{}
\vspace{.15cm}
\\
{}+
\omega(d(x_1,x_2)) \min\{\sqrt{|\xi_1|^2+|\eta_1|^2}, \sqrt{|\xi_2|^2+|\eta_2|^2}\}
\end{array}
\end{equation}
for every $x_1,x_2\in X$, $\xi_1,\xi_2\in\Xi$, $\eta_1,\eta_2\in\R$.

\par
By definition we have 
\begin{equation}\label{500.5}
\langle f,\delta_{u_n}\rangle = \int_{X^{n,a}}f(x,u_n,1)\,d\lambda + 
\int_{X^{n,s}}f(x,u_n,1)\,d\lambda\,. 
\end{equation}
By \eqref{pndef1} the first integral in the right-hand side can be written as
\begin{equation}\label{501}
\begin{array}{c}
\displaystyle
\int_{X^{n,a}}f(x,u_n,1)\,d\lambda= \sum_{ij}\int_{A^{n,Y}_{ij}}f(x,\xi^{n,Y}_j\!,1)\, 
d\lambda +
\sum_{ij}\int_{A^{n,\infty}_{ij}}f(x,\xi^{n,\infty}_j/\e^n_i,1)\, d\lambda =
\\
\displaystyle
=\sum_{ij} f(x^{n,a}_i,\xi^{n,Y}_j\!,1)\lambda(A^{n,Y}_{ij}\!) + 
\sum_{ij}f(x^{n,a}_i,\xi^{n,\infty}_j/\e^n_i,1)\lambda(A^{n,\infty}_{ij}) +r^{a,1}_n,
\end{array}
\end{equation}
where $x^{n,a}_i$ are arbitrary points in $A^{n,a}_i$ and the remainder $r^{a,1}_n$ 
tends to $0$ as a consequence of \eqref{diam2}, \eqref{choice1}, \eqref{choice2}, and 
\eqref{500}, which lead to the estimate
\begin{eqnarray*}
&\displaystyle
|r^{a,1}_n|\le \omega(\sigma_n) \sum_{ij}\sqrt{1+|\xi^{n,Y}_j\!|^2} \, 
\lambda(A^{n,Y}_{ij}\!) +
\omega(\sigma_n) \sum_{ij}\frac{1}{\e^n_i}\sqrt{(\e^n_i)^2+|\xi^{n,\infty}_j|^2}\,
\lambda(A^{n,\infty}_{ij})\le
\\
&\displaystyle
\le \sigma_n\, \omega(\sigma_n)\, \lambda(X^{n,a})+
 \omega(\sigma_n)\, \lambda^Y\!(X^{n,a})+
\omega(\sigma_n)\, \lambda^{\infty}(X^{n,a})
\sqrt{1+\sigma_n^2}\,.
\end{eqnarray*}
On the other hand by \eqref{choice1} and \eqref{choice2} we have
\begin{equation}\label{502}
\begin{array}{c}
\displaystyle
\sum_{ij}f(x^{n,a}_i,\xi^{n,Y}_j\!,1)\,\lambda(A^{n,Y}_{ij}\!)
+\sum_{ij}f(x^{n,a}_i,\xi^{n,\infty}_j/\e^n_i,1)\,\lambda(A^{n,\infty}_{ij})=
\\
\displaystyle
=\sum_{ij}(1-\delta^n_i)\, f(x^{n,a}_i,\xi^{n,Y}_j\!,1) 
\,\mu^Y\!(A^{n,a}_i{\times}B^{n,Y}_j\!)+{}
\\
\displaystyle
+\sum_{ij}f(x^{n,a}_i,{\xi^{n,\infty}_j},\e^n_i)\, 
\mu^\infty(A^{n,a}_i{\times}B^{n,\infty}_j)=
\\
\displaystyle
=\int_{X{\times}\Xi} f(x,\xi,1)\,d\mu^Y\!(x,\xi) +\int_{X^a{\times}\Sigma_\Xi} 
f(x,\xi,0)\,d\mu^\infty(x,\xi)+ r^{a,2}_n\,,
\end{array}
\end{equation}
where the remainder $r^{a,2}_n$ tends to $0$ as a consequence of \eqref{diam1}, 
\eqref{702}, \eqref{diam2}, and \eqref{500}, which lead to the estimate
$$
|r^{a,2}_n|\le \big(2a\sigma_n+\omega(\sigma_n)\big) \big( \lambda^Y\! (X) +
\lambda^{\infty}(X^a) \sqrt{1+\sigma_n^2}\big)+2a \sigma_n\,.
$$
{}From \eqref{501} and  \eqref{502} we obtain
\begin{equation}\label{503}
\int_{X^{n,a}}f(x,u_n,1)\,d\lambda=
\int_{X{\times}\Xi} f(x,\xi,1)\,d\mu^Y\!(x,\xi) +\int_{X^a{\times}\Sigma_\Xi} 
f(x,\xi,0)\,d\mu^\infty(x,\xi) +r^{a}_n\,,
\end{equation}
where $r^{a}_n:=r^{a,1}_n+r^{a,2}_n$ tends to $0$.

\par
By \eqref{pndef4} and \eqref{defcnk} the second integral in the right-hand side of 
\eqref{500.5} can be written as
\begin{equation}\label{507}
\displaystyle
\int_{X^{n,s}}f(x,u_n,1)\,d\lambda= \sum_{ij} 
f(x^{n,s}_i,c^{n}_{ij}\xi^{n,\infty}_j,1)\lambda(A^{n,s}_{ij})+r^{s,1}_n\, ,
\end{equation}
where $x^{n,s}_i$ are arbitrary points in $X^{n,s}_i$ and the remainder $r^{s,1}_n$ 
tends to zero as a consequence of \eqref{diam-2}, \eqref{diam0}, \eqref{defcnk}, 
\eqref{choice0.5}, and \eqref{500}, which lead to the estimate
$$
|r^{s,1}_n|\le \omega(\sigma_n)\, \lambda^{\infty}(X^s)\sqrt{1+\sigma_n^2}\,.
$$
On the other hand by \eqref{defcnk} 
\begin{equation}\label{508}
\begin{array}{c}
\displaystyle
\sum_{ij} f(x^{n,s}_i,c^{n}_{ij}\xi^{n,\infty}_j,1)\lambda(A^{n,s}_{ij})=\sum_{ij} 
f(x^{n,s}_i,\xi^{n,\infty}_j,1/c^n_{ij})\,\mu^\infty(X_i^{n,s}{\times}B^{n,\infty}_j)=
\\
\displaystyle
=\int_{X^{s}{\times}\Sigma_\Xi} f(x,\xi,0)\,d\mu^\infty(x,\xi)+ r^{s, 2}_n
\,,
\end{array}
\end{equation}
where the remainder $r^{s,2}_n$ tends to $0$ as a consequence of \eqref{diam1}, 
\eqref{diam-2}, \eqref{choice0.5}, and \eqref{500}, which lead to the estimate
$$
|r^{s,2}_n|\le \big(a\sigma_n +\omega(\sigma_n)\sqrt{1+\sigma_n^2} \big)
\lambda^{\infty}(X^s)\,.
$$
{}From \eqref{507} and  \eqref{508} we obtain
\begin{equation}\label{509}
\int_{X^{n,s}}f(x,u_n,1)\,d\lambda=
\int_{X^{s}{\times}\Sigma_\Xi} f(x,\xi,0)\,d\mu^\infty(x,\xi)+ r^{s}_n\,,
\end{equation}
where $r^{s}_n:=r^{s,1}_n+r^{s,2}_n$ tends to $0$.

\par
{}From \eqref{mu1muinf}, \eqref{500.5}, \eqref{503}, and \eqref{509} we obtain 
\eqref{499}, which concludes the proof of the theorem.
\end{proof}

\par
\begin{remark}\label{normpn}
If $u_n$ is a sequence in $L^1(X;\Xi)$ such that $\delta_{u_n}\wto \mu$ weakly$^*$ 
in $GY(X;\Xi)$, then
$$
\int_X\sqrt{1+|u_n|^2}\, d\lambda \longrightarrow \|\mu\|_* 
$$
by Lemma~\ref{lemma:mass} and Remark~\ref{convnorms}.
\end{remark}

\end{section}

\begin{section}{The notion of barycentre}

In this section we study some properties of the barycentre of a generalized Young 
measure.

\par
\begin{definition}
The {\em barycentre\/} of a generalized Young measure $\mu\in GY(X;\Xi)$ is the 
measure $\bary(\mu)\in M_b(X;\Xi)$ defined by
$$
\bary(\mu) =\pi_X(\xi\,\mu)\,.
$$
\end{definition}

\par
\begin{remark}\label{bar9}
By Definition~\ref{def:pihmu} a measure $p\in M_b(X;\Xi)$ coincides with 
$\bary(\mu)$ if and only if
\begin{equation}\label{2bary}
\int_X \varphi{\,\cdot\,} dp=\langle\varphi(x) {\,\cdot\,}\xi,
\mu(x,\xi,\eta)\rangle
\end{equation}
for every $\varphi\in C(X;\Xi)$. By approximation we can prove that the same equality 
holds for every bounded Borel function $\varphi:X\to\Xi$.
\end{remark}

\par
\begin{remark}\label{bar10}
Let $\mu=\ol\mu^Y\! + \hat\mu^\infty$ be the decomposition of 
Theorem~\ref{thm:decomp}, let $(\mu^{x,Y}\!)_{x\in X}$ and 
$(\mu^{x,\infty})_{x\in X}$ be the families of probability measures introduced in 
Remark~\ref{disint}, and let $\lambda^\infty:=\pi_X(\mu^\infty)$. 
We consider the functions $u^Y\!\colon\Om\to\Xi$ and $u^\infty\colon\Om\to\Xi$
defined by
$$
u^Y\!(x):=\int_\Xi \xi\,d\mu^{x,Y}\!(\xi) \qquad\hbox{and}  \qquad u^\infty(x):=
\int_{\Sigma_\Xi} \xi\,d\mu^{x,\infty}(\xi)\,.
$$
Then $\bary(\mu)=u^Y\! + u^\infty\,\lambda^\infty$.
In particular, if $\mu=\ol\mu^Y$, then $\bary(\mu)=u^Y\!\in L^1(X;\Xi)$.
Therefore, $\bary(\delta_u)=u$ for every $u\in L^1(X;\Xi)$.
It follows immediately from Definition~\ref{GYMp} and \eqref{2bary} that we have also
$\bary(\delta_p)=p$ for every $p\in M_b(X;\Xi)$.
\end{remark}

\par
\begin{remark}\label{bar11}
{}From Remark~\ref{piXhmu} we obtain
$$
\|\bary(\mu)\|\le \|\mu\|_*\,.
$$
If $\mu_k\wto\mu$ weakly$^*$ in $GY(X;\Xi)$, then $\bary(\mu_k)\wto \bary(\mu)$ 
weakly$^*$ in $M_b(X;\Xi)$.
If $\mu=\ol\mu^Y$ with $\mu^Y\in Y^r(X;\Xi)$ for some $r>1$, then 
Remark~\ref{bar10} implies that $\bary(\mu)\in L^r(X;\Xi)$ and
$$
\|\bary(\mu)\|_r\le \Big( \int_{X{\times}\Xi} |\xi|^r\,d\mu^Y\!(x,\xi)\Big)^{1/r}= 
 \langle\{P_r\},\mu\rangle^{1/r}\,,
$$
where $\|\cdot\|_r$ denotes the norm in $L^r(X;\Xi)$ and $\{P_r\}$ is the homogeneous 
function defined in~\eqref{mom1}.
\end{remark}

\par
We now prove the {\em Jensen inequality\/} for generalized Young measures.

\par
\begin{theorem}\label{Jens}
Let $f\colon X{\times}\Xi{\times}\R\to \R\cup\{+\infty\}$ be a Borel function such that
$(\xi,\eta)\mapsto f(x,\xi,\eta)$ is positively one-homogeneous, convex, and lower 
semicontinuous for every $x\in X$ and satisfies the inequality
$$
f(x,\xi,\eta)\ge -c\, \sqrt{|\xi|^2+|\eta|^2}
$$
for some constant $c$.
Then
\begin{equation}\label{jensen}
\langle f, \delta_{\bary(\mu)}\rangle \le \langle f, \mu\rangle
\end{equation}
for every $\mu\in GY(X;\Xi)$.
\end{theorem}

\par
\begin{proof}
Let us fix $\mu\in  GY(X;\Xi)$, let $p:=\bary(\mu)$, and let 
$\sigma\in M_b^+(X)$ be such that $p<<\sigma$ and $\lambda<<\sigma$. 
We consider an increasing sequence of functions $f_k$ converging to $f$ such that each
$f_k$ has the form 
$$
f_k(x,\xi,\eta)=\sup_{1\le i\le k}\{ a_i(x){\,\cdot\,}\xi+b_i(x)\eta \}
$$
with $a_i\colon X\to \Xi$ and $b_i\colon X\to\R$ bounded $\sigma$-measurable 
functions (see, e.g., \cite[Theorem~2.2.4]{But}). For every $k$ there exists a Borel 
partition $(B^k_i)_{1\le i\le k}$ such that
$$
\int_{X} f_k(x,{\textstyle  \frac{dp}{d\sigma},\frac{d\lambda}{d\sigma}})\, d\sigma 
=\sum_{i=1}^k \Big\{
\int_{B^k_i} a_i{\,\cdot\,}dp+ \int_{B^k_i} b_i\, d\lambda \Big\}\,.
$$
By (\ref{muproj}) and (\ref{2bary}) we obtain
\begin{eqnarray*}
& \displaystyle
 \int_{B^k_i} a_i{\,\cdot\,} dp+ \int_{B^k_i} b_i\, d\lambda=
\langle (a_i(x){\,\cdot\,}\xi+b_i(x)\eta)1_{B^k_i}(x), \mu(x,\xi,\eta) \rangle \le
\\
& \le \langle f(x,\xi,\eta)1_{B^k_i}(x), \mu(x,\xi,\eta) \rangle \,.
\end{eqnarray*}
Summing over $i$ we get
$$
\int_{X} f_k(x, {\textstyle\frac{dp}{d\sigma}, \frac{d\lambda}{d\sigma}})\, d\sigma \le
\langle f(x,\xi,\eta), \mu(x,\xi,\eta) \rangle\,,
$$
and taking the limit with respect to $k$ gives inequality~(\ref{jensen}).
\end{proof}

\par
\begin{remark}\label{convexification}
Let $f\colon X{\times}\Xi{\times}\R\to [0,+\infty]$ be a Borel function such that 
$(\xi,\eta)\mapsto f(x,\xi,\eta)$ is positively one-homogeneous for every $x\in X$, and let 
${\rm \ol{co}}\, f$ be the lower semicontinuous convex envelope of $f$ with respect to 
$(\xi,\eta)$. By applying \eqref{jensen} to ${\rm \ol{co}}\, f$ we obtain
$$
\langle {\rm \ol{co}}\, f, \delta_{\bary(\mu)}\rangle \le \langle f, \mu\rangle
$$
for every $\mu\in GY(X;\Xi)$.
\end{remark}

\par
The opposite inequality requires special conditions on $f$ and $\mu$, as shown in the 
following lemma, that will be used in~\cite{DM-DeS-Mor-Mor}.

\par
\begin{lemma}\label{lm:supporto}
Let $\mu\in GY(X;\Xi)$, let  $f\colon X{\times}\Xi{\times}\R\to [0,+\infty]$ be a Borel 
function such that $(\xi,\eta)\mapsto f(x,\xi,\eta)$ is positively one-homogeneous for every 
$x\in X$, and let ${\rm \ol{co}}\, f$ be the lower semicontinuous convex envelope of $f$ 
with respect to $(\xi,\eta)$. Assume that 
$ \langle f,\mu\rangle\le \langle{\rm \ol{co}}\, f, \delta_{\bary(\mu)}\rangle <+\infty$.
Then $\supp\,\mu$ is contained in the closure of $\{ f={\rm \ol{co}}\, f\}$.
\end{lemma}

\par
\begin{proof}
Using the hypothesis and \eqref{jensen} we obtain
$$
\langle {\rm \ol{co}}\, f, \delta_{\bary(\mu)}\rangle \le \langle {\rm \ol{co}}\, f, 
\mu\rangle \le \langle f,\mu\rangle \le \langle {\rm \ol{co}}\, f, \delta_{\bary(\mu)}\rangle\,,
$$
hence, $\langle f-{\rm \ol{co}}\, f,\mu\rangle=0$.
Since $f-{\rm \ol{co}}\, f$ and $\mu$ are nonnegative, we conclude that
$\supp\,\mu$ is contained in the closure of $\{ f={\rm \ol{co}}\, f\}$.
\end{proof}

\end{section}

\begin{section}{Compatible systems of generalized Young measures}

Let $A\subset\R$ and let $t\mapsto p(t)$ be a function from $A$ into $M_b(X ;\Xi)$. 
For every finite sequence $t_1 < t_2 < \dots < t_m$ in $A$ we consider the measure 
$(p(t_1),\dots,p(t_m))\in M_b(X ;\Xi^m)$ and the corresponding generalized Young 
measure 
\begin{equation}\label{delta1m}
\mu_{t_1\dots t_m}:=\delta_{(p(t_1),\dots,p(t_m))}\in GY(X ;\Xi^m)
\end{equation}
introduced in Definition~\ref{GYMp}, with $\Xi$ replaced by $\Xi^m$. To describe an 
important property of this family of generalized Young measures it is convenient to 
introduce the following definition.

\par
\begin{definition}\label{pits}
If $\{s_1,s_2,\dots,s_n\}\subset\{t_1,t_2,\dots,t_m\}\subset \R$, with 
$s_1<s_2<\dots<s_n$ and $t_1 < t_2 < \dots < t_m$, we define the projection 
$\pi_{s_1\dots s_n}^{t_1\dots t_m}\colon X{\times}\Xi^m{\times}\R\to 
X{\times}\Xi^n{\times}\R$ by 
$$
\pi_{s_1\dots s_n}^{t_1\dots 
t_m}(x,\xi_{t_1},\dots,\xi_{t_m},\eta)=(x,\xi_{s_1},\dots,\xi_{s_n},\eta)\,.
$$
\end{definition}

\par
\begin{remark}\label{rem:proj}
It is easy to see that the family of generalized Young measures \eqref{delta1m} satisfies 
the {\it compatibility condition\/}
\begin{equation}\label{compatib}
\mu_{s_1\dots s_n}=\pi_{s_1\dots s_n}^{t_1\dots t_m}(\mu_{t_1\dots t_m})
\end{equation}
whenever $\{s_1,s_2,\dots,s_n\}$ and $\{t_1,t_2,\dots,t_m\}$ are as in 
Definition~\ref{pits}.
\end{remark}

\par
This motivates the following definition.

\par
\begin{definition}\label{defSGY}
A {\it compatible system of generalized Young measures\/} on $X$, with values in a
finite dimensional Hilbert space $\Xi$ and with time set $A\subset \R$, is a family
$\mu=(\mu_{t_1\dots t_m})$ of generalized Young measures 
$\mu_{t_1\dots t_m}\in GY(X ;\Xi^m)$ satisfying the compatibility condition
(\ref{compatib}), with $t_1, \dots,t_m$ running over all finite sequences of 
elements of $A$ with $t_1 < t_2 < \dots < t_m$. The space of all such systems
is denoted by $SGY(A,X ;\Xi)$ and is equipped with the weakest topology such 
that the maps $\mu\mapsto \mu_{t_1\dots t_m}$ from $SGY(A,X ;\Xi)$ into 
$GY(X;\Xi^m)$, endowed with the weak$^*$ topology, are continuous for every 
$m$ and every finite sequence $t_1, \dots,t_m$ in $A$ with $t_1 < t_2 < \dots < t_m$.
Although this topology is not induced by duality, we shall refer to it as the 
{\em weak$^*$ topology\/} of $SGY(A,X ;\Xi)$.
\end{definition}

\par
\begin{definition}\label{defSGYp}
Given a function $t\mapsto p(t)$ from $A\subset \R$ into $M_b(X ;\Xi)$, the family 
$(\mu_{t_1\dots t_m})$ defined in (\ref{delta1m}) is called the {\em compatible system 
of generalized Young measures associated with\/} $t\mapsto p(t)$.
\end{definition}

\par
The compatibility condition \eqref{compatib} implies that the barycentre of 
$\mu_{t_1\dots t_m}$ is completely determined by the barycentres of 
$\mu_{t_1},\dots,\mu_{t_m}$.

\par
\begin{proposition}\label{rembary}
Let $\mu\in SGY(A,X;\Xi)$. Then
$$
\bary(\mu_{t_1\dots t_m})=(\bary(\mu_{t_1}),\dots,\bary(\mu_{t_m}))
$$
for every finite sequence $t_1, \dots,t_m$ in $A$ with $t_1 < t_2 < {\dots < t_m}$.
\end{proposition}

\par
\begin{proof}
Let $(p_1,\dots,p_m):=\bary(\mu_{t_1\dots t_m})$ and $q_i=\bary(\mu_{t_i})$ for 
$i=1, \dots,m$.
Using \eqref{2bary} for every $(\varphi_1,\dots,\varphi_m)\in C(X;\Xi^m)$ we have
\begin{eqnarray}
&\displaystyle
\sum_{i=1}^m \int_X \varphi_i
{\,\cdot\,}dp_i=\sum_{i=1}^m 
\langle \varphi_i
(x){\,\cdot\,} \xi_i,\mu_{t_1\dots t_m}(x, \xi_1,\dots,\xi_m,\eta)\rangle\,,
\label{901}
\\
&\displaystyle
\int_X \varphi_i
{\,\cdot\,}dq_i=\langle \varphi_i
(x){\,\cdot\,} \xi_i,\mu_{t_i}(x, \xi_i,\eta)\rangle\qquad\hbox{for } i=1, \dots,m
\,.\label{902}
\end{eqnarray}
The compatibility condition \eqref{compatib} implies that
$$
\langle \varphi_i
(x){\,\cdot\,} \xi_i,\mu_{t_1\dots t_m}(x, \xi_1,\dots,\xi_m,\eta)\rangle
=
\langle \varphi_i(x){\,\cdot\,} \xi_i,\mu_{t_i}(x, \xi_i,\eta)\rangle\,,
$$
hence \eqref{901} and \eqref{902} yield
$$
\sum_{i=1}^m \int_X \varphi_i{\,\cdot\,}dp_i=
\sum_{i=1}^m \int_X \varphi_i{\,\cdot\,}dq_i
$$
for every $(\varphi_1,\dots,\varphi_m)\in C(X;\Xi^m)$.
This gives $p_i=q_i$ for $i=1, \dots,m$.
\end{proof}

\par
The notion of left continuity, introduced in the next definition, is very useful in the 
applications.

\par
\begin{definition}\label{defSGYlc}
A system $\mu\in SGY(A,X ;\Xi)$ is said to be {\it left continuous\/} if for every finite 
sequence $t_1,\dots,t_m$ in $A$ with $t_1 <\dots < t_m$ the following  continuity 
property holds:
\begin{equation} \label{leftlim}
 \mu_{s_1\dots s_m}\wto \mu_{t_1\dots t_m} \quad\hbox{weakly}^*\hbox{ in }
 GY(X ;\Xi^m)
\end{equation}
as $s_i\to t_i$, with $s_i\in A$ and $s_i\le t_i$.
\end{definition}

\par
The following theorem proves the weak$^*$ compactness of the subsets of 
$SGY(A,X;\Xi)$ defined by imposing  bounds on the norms of $\mu_t$ for every
$t\in A$.

\par
\begin{theorem}\label{compactness10}
For every function $C\colon A\to[0,+\infty)$ the set
\begin{equation}\label{set-cpt}
\{ \mu\in SGY(A,X;\Xi): \ \|\mu_t\|_*\le C(t) \ \text{for every } t\in A\}
\end{equation}
is weakly$^*$ compact in $SGY(A,X ;\Xi)$.
\end{theorem}

\par
To prove the theorem we need the following lemma which provides an estimate of the 
norm $\|\mu_{t_1\dots t_m}\|_*$ in terms of the norms $ \|\mu_{t_i}\|_*$.

\par
\begin{lemma}\label{900}
For every $\mu\in SGY(A,X ;\Xi)$ we have
$$
\|\mu_{t_1\dots t_m}\|_*   \le \sum_{i=1}^m \|\mu_{t_i}\|_*
$$
for every finite sequence $t_1, \dots,t_m$ in $A$ with $t_1 < t_2 < \dots < t_m$. 
\end{lemma}

\par
\begin{proof}
By Remark~\ref{M+} and by the compatibility condition (\ref{compatib})  we have
\begin{eqnarray*}
& \|\mu_{t_1\dots t_m}\|_*=
\langle |(\xi_1,\dots,\xi_m,\eta)|,
\mu_{t_1\dots t_m}(x,\xi_1,\dots,\xi_m,\eta)\rangle \le 
\\
&\displaystyle
\le  \sum_{i=1}^m 
\langle |(\xi_i,\eta)|,
\mu_{t_1\dots t_m}(x,\xi_1,\dots,\xi_m,\eta)\rangle =
 \sum_{i=1}^m 
\langle |(\xi_i,\eta)|,
\mu_{t_i}(x,\xi_i,\eta)\rangle=
\sum_{i=1}^m \|\mu_{t_i}\|_* \,,
\end{eqnarray*}
which concludes the proof.
\end{proof}

\par
\begin{proof}[Proof of Theorem~\ref{compactness10}]
By Lemma~\ref{900} for every function $C\colon A\to [0,+\infty)$ the set defined in
\eqref{set-cpt} is contained in the set of all $\mu\in SGY(A,X;\Xi)$ such that
$$
\|\mu_{t_1\dots t_m}\|_*\le \sum_{i=1}^m C(t_i)
$$
for every finite sequence $t_1, \dots,t_m$ in $A$ with $t_1 < t_2 < \dots < t_m$. As the 
topology in $SGY(A,X;\Xi)$ is induced by the product of the weak$^*$ topologies of the 
spaces $GY(X;\Xi^m)$ corresponding to the projections $\mu_{t_1\dots t_m}$, the set 
(\ref{set-cpt}) is compact in the weak$^*$ topology of $SGY(A,X;\Xi)$ by 
Tychonoff's Theorem.
\end{proof}

\par
\begin{remark}\label{Afinite}
If $A=\{a_0,a_1,\dots,a_k\}$, with $a_0<a_1<\dots<a_k$, then for every 
$\mu\in GY(X;\Xi^{k+1})$ there exists a unique system $\mu^A\in SGY(A,X ;\Xi)$
such that $\mu^A_{a_0\dots a_k}=\mu$. This system is defined by 
$$
\mu^A_{t_1\dots t_m}=\pi^{a_0\dots a_k}_{t_1\dots t_m}(\mu)
$$
for every $\{t_1,t_2,\dots,t_m\}\subset \{a_0,a_1,\dots,a_k\}$ with
$t_1 < t_2 < \dots < t_m$.
\end{remark}

\par
The notion of  piecewise constant interpolation will be useful in the application to 
evolution problems.

\par
\begin{definition}\label{interpol}
Let $A=\{a_0,a_1,\dots,a_k\}$, with $a_0<a_1<\dots<a_k$. For every $t_1, \dots,t_m$ in 
$[a_0,a_k]$ with $t_1 < t_2 < \dots < t_m$ let $\rho_{t_1\dots t_m}\colon 
X{\times}\Xi^{k+1}{\times}\R \to X{\times}\Xi^m{\times}\R$ be defined by 
$$
\rho_{t_1\dots t_m}(x,\xi_{a_0},\dots,\xi_{a_k},\eta):=(x,\xi_{t_1},\dots,\xi_{t_m},\eta)\,,
$$
with $\xi_{t_i}=\xi_{a_j}$, where $j$ is the largest index such that $a_j\le t_i$.
For every $\mu\in GY(X ;\Xi^{k+1})$ the {\it piecewise constant interpolation\/} 
$\mu^{[A]}$ of $\mu$ is the element of $SGY([a_0,a_k],X ;\Xi)$ defined by
\begin{equation}\label{piecconst}
\mu^{[A]}_{t_1\dots t_m}:=\rho_{t_1\dots t_m}(\mu)
\end{equation}
for every $t_1, \dots,t_m$ in $[a_0,a_k]$ with $t_1 < t_2 < \dots < t_m$.
\end{definition}

\par
\begin{remark}\label{rem:interpol}
It is easy to check that $\rho_{s_1\dots s_n}=\pi_{s_1\dots s_n}^{t_1\dots 
t_m}\circ\rho_{t_1\dots t_m}$ whenever 
$\{s_1,s_2,\dots,s_n\}\!\subset\!\{t_1,t_2,\dots,t_m\}\subset [a_0,a_k]$, with 
$s_1<s_2<\dots<s_n$ and $t_1 < t_2 < \dots < t_m$.
Therefore the family of generalized Young measures
$(\mu^{[A]}_{t_1\dots t_m})$ defined by \eqref{piecconst} satisfies the compatibility 
condition (\ref{compatib}).
\end{remark}

\end{section}

\begin{section}{The notion of variation}

In this section we study the notion of variation on a time interval  of a compatible system 
of generalized Young measures, and prove a compactness theorem which extends Helly's 
Theorem.

\par
\begin{definition}\label{variation}
Given a set $A\subset\R$, the  {\it variation\/} of $\mu\in SGY(A,X ;\Xi)$ on the time 
interval $[a,b]$, with $a$, $b\in A$, is defined as 
$$
{\rm Var}(\mu;a,b):= \sup \sum_{i=1}^k 
\langle|\xi_i-\xi_{i\!-\!1}|, \mu_{t_0t_1\dots t_k}(x,\xi_0,\dots,\xi_k, \eta) \rangle\,,
$$
where the supremum is taken over all finite families $t_0,t_1,\dots,t_k$ in $A$ such that 
$a=t_0<t_1<\dots<t_k=b$ (with the convention ${\rm Var}(\mu;a,b)=0$ if $a=b$).
\end{definition}

\par
\begin{remark}\label{rem:344}
If $\mu$ is the compatible family of generalized Young measures associated with a
function $t\mapsto p(t)$ from $A$ into $M_b(X ;\Xi)$ according to \eqref{delta1m}, 
then
${\rm Var}(\mu;a,b)$ reduces to the variation of $t\mapsto p(t)$ on ${[a,b]\cap A}$. 
\end{remark}

\par
\begin{remark}\label{rem:345}
Returning to the general case, the compatibility condition (\ref{compatib}) yields
$$
{\rm Var}(\mu;a,b) = \sup 
\sum_{i=1}^k \langle|\xi_i-\xi_{i\!-\!1}|, \mu_{t_{i\!-\!1}t_i}(x,\xi_{i\!-
\!1},\xi_i,\eta)\rangle\,,
$$
where the supremum is taken over all finite families $t_0,t_1,\dots,t_k$ in $A$ such that 
$a=t_0<t_1<\dots<t_k=b$.
\end{remark}

\par
\begin{remark}\label{rem:346}
If $t_1,t_2,t_3\in A$ and $t_1<t_2<t_3$, by the compatibility condition (\ref{compatib}) 
and by the triangle inequality we have
\begin{eqnarray*}
& \langle|\xi_3 -\xi_1| , \mu_{t_1t_3}(x,\xi_1,\xi_3,\eta) \rangle= 
 \langle|\xi_3 -\xi_1| ,  \mu_{t_1t_2t_3}(x,\xi_1,\xi_2,\xi_3,\eta) \rangle \le 
 \\
& \le   \langle|\xi_3 -\xi_2| ,  \mu_{t_1t_2t_3}(x,\xi_1,\xi_2,\xi_3,\eta) \rangle +
 \langle|\xi_2 -\xi_1| ,  \mu_{t_1t_2t_3}(x,\xi_1,\xi_2,\xi_3,\eta) \rangle =
 \\
 & = \langle|\xi_3 -\xi_2| ,  \mu_{t_2t_3}(x,\xi_2,\xi_3,\eta) \rangle +
 \langle|\xi_2 -\xi_1| ,  \mu_{t_1t_2}(x,\xi_1,\xi_2,\eta) \rangle\,.
\end{eqnarray*} 
Using this inequality it is easy to deduce from Remark~\ref{rem:345} that
\begin{equation}\label{additive}
{\rm Var}(\mu;a,c)= {\rm Var}(\mu;a,b) + {\rm Var}(\mu;b,c) 
\end{equation}
for every $a, b, c\in A$ with $a \le b \le c$. This implies in particular that the function 
$t\mapsto {\rm Var}(\mu;a,t)$ is nondecreasing on $A\cap[a,+\infty)$.
\end{remark}

\par
\begin{remark}\label{rem:347}
If $A=\{a_0,\dots,a_k\}\subset\R$ is a finite set,  with $a_0<a_1<\dots <a_k$, 
$\mu\in GY(X;\Xi^{k+1})$, and $\mu^A\in SGY(A,X ;\Xi)$ is the associated 
system defined in Remark~\ref{Afinite}, it follows from (\ref{additive}) that 
$$
\begin{array}{c}
\displaystyle
 {\rm Var}(\mu^A;a_0,a_k)=\sum_{i=1}^k 
 \langle|\xi_i-\xi_{i\!-\!1}|, \mu^A_{a_{i\!-\!1}a_i}(x,\xi_{i\!-\!1},\xi_i,\eta)\rangle=
 \\
 \displaystyle
 =\sum_{i=1}^k 
 \langle|\xi_i-\xi_{i\!-\!1}|, \mu(x,\xi_{0},\dots,\xi_k,\eta)\rangle
 \,.
 \end{array}
$$
It is easy to see that, if  $\mu^{[A]}\in SGY([a_0,a_k],X ;\Xi)$ is the piecewise constant 
interpolation of $\mu$ defined by (\ref{piecconst}), then
$$
\begin{array}{c}
\displaystyle
{\rm Var}(\mu^{[A]};a_0,a_k)={\rm Var}(\mu^A;a_0,a_k)=
\sum_{i=1}^k 
 \langle|\xi_i-\xi_{i\!-\!1}|, \mu^A_{a_{i\!-\!1}a_i}(x,\xi_{i\!-\!1},\xi_i,\eta)\rangle
 \\
 \displaystyle
 =\sum_{i=1}^k 
 \langle|\xi_i-\xi_{i\!-\!1}|, \mu(x,\xi_{0},\dots,\xi_k,\eta)\rangle
 \,.
 \end{array}
$$
\end{remark}

\par
\begin{definition}\label{variationh}
Let $h\colon\Xi\to[0,+\infty)$ be a positively one-homogeneous function satisfying the 
triangle inequality.
Given a set $A\subset\R$, the  {\it $h$-variation\/} of $\mu\in SGY(A,X ;\Xi)$ on the 
time interval $[a,b]$, with $a$, $b\in A$, is defined as 
$$
{\rm Var}_h(\mu;a,b):= \sup \sum_{i=1}^k 
\langle h(\xi_i-\xi_{i\!-\!1}), \mu_{t_0t_1\dots t_k}(x,\xi_0,\dots,\xi_k, \eta) \rangle\,,
$$
where the supremum is taken over all finite families $t_0,t_1,\dots,t_k$ in $A$ such that 
$a=t_0<t_1<\dots<t_k=b$ (with the convention ${\rm Var}_h(\mu;a,b)=0$ if $a=b$).
\end{definition}

\par
\begin{remark}\label{rem-Varh}
It is well known that every positively one-homogeneous function 
$h\colon\Xi\to[0,+\infty)$
satisfying the triangle inequality is continuous and satisfies an estimate of the form
$h(\xi)\le c\,|\xi|$ for some constant $c$.
It follows that ${\rm Var}_h(\mu;a,b)\le c\,{\rm Var}(\mu;a,b)$.
It is easy to see that all properties of ${\rm Var}(\mu;a,b)$ proved so far can be 
extended to ${\rm Var}_h(\mu;a,b)$.
\end{remark}

\par
Using the compatibility condition it is easy to prove the following lemma.

\par
\begin{lemma}
Let $T>0$ and let $\mu\in SGY([0,T], X;\Xi)$ with ${\rm Var}(\mu;0,T)<+\infty$. For 
every $f\in C^\hom_L(X{\times}\Xi{\times}\R)$ the function $t\mapsto \langle 
f,\mu_t\rangle$ has bounded variation on $[0,T]$.
\end{lemma}

\par
The proof is omitted, since it is similar to the proof of the following lemma, which will be 
used in Theorem~\ref{BV-weakder}.

\par
\begin{lemma}\label{lm:bv-mut}
Let $T>c>0$ and let $\mu\in SGY([0,T], X;\Xi)$ with ${\rm Var}(\mu;0,T)<+\infty$. 
For every $f\in C^\hom_L(X{\times}\Xi^2{\times}\R)$ the function 
$\Phi_c^f(t):= \langle f,\mu_{t, t+c}\rangle$ has bounded variation on $[0,T-c]$.
\end{lemma}

\par
\begin{proof}
Let $V(t):={\rm Var}(\mu;0,t)$ for every $t\in[0,T]$.
Let us fix $f\in C^\hom_L(X{\times}\Xi^2{\times}\R)$ and let $a$ be a constant 
satisfying~\eqref{Lip}. Let $t_1,t_2$ with $0\le t_1<t_1+c<t_2<t_2+c\le T$. 
Using the compatibility condition \eqref{compatib} and \eqref{additive}, we obtain
$$
\begin{array}{c}
|\Phi_c^f(t_2)-\Phi_c^f(t_1)| =
| \langle f(x,\xi_2,\xi_2',\eta) -  f(x,\xi_1,\xi_1',\eta) , \mu_{t_1, t_1+c, t_2, 
t_2+c}(x,\xi_1,\xi_1',\xi_2,\xi_2',\eta)\rangle | \le
\smallskip
\\
\le 
a \langle |\xi_2-\xi_1|+ |\xi_2'-\xi_1'| , \mu_{t_1, t_1+c, t_2, 
t_2+c}(x,\xi_1,\xi_1',\xi_2,\xi_2',\eta)\rangle =
\smallskip
\\
= a \langle |\xi_2-\xi_1| , \mu_{t_1, t_2}(x,\xi_1,\xi_2,\eta)\rangle
+ a \langle |\xi_2'-\xi_1'| , \mu_{ t_1+c, t_2+c}(x,\xi_1',\xi_2',\eta)\rangle
\le
\smallskip
\\
\le V(t_2)-V(t_1)+V(t_2+c)-V(t_1+c)\,.
\end{array}
$$
The same inequality can be proved if $0\le t_1<t_2\le t_1+c<t_2+c\le T$.
As $V$ is nondecreasing, we conclude that the total variation of $\Phi_c^f$ on $[0,T-c]$ 
is less than or equal to $V(T-c)+V(T)$.
\end{proof}

\par
The following result can be considered as a version of Helly's Theorem for compatible 
systems of generalized Young measures. Note that this is a sequential compactness result, 
in contrast with Theorem~\ref{compactness10}.

\par
\begin{theorem}\label{Helly} 
Let $T>0$ and let $\mu^k$ be a sequence in $SGY([0,T],X ;\Xi)$  such that
\begin{eqnarray}
&\displaystyle \sup_{k} {\rm Var}(\mu^k;0,T)\le C\,,\label{boundvar}
\\
&\displaystyle \sup_{k}\|\mu^k_{t_0}\|_*\le C_*  \,,\label{boundnorm}
\end{eqnarray}
for some $t_0\in [0,T]$ and some finite constants $C$ and $C_*$.
Then there exist a subsequence, still denoted $\mu^k$, a set $\Theta\subset [0,T]$, 
containing $0$ and with $[0,T]\setmeno \Theta$ at most countable, and a left continuous 
$\mu\in SGY([0,T],X ;\Xi)$, with
\begin{eqnarray}
& {\rm Var}(\mu;0,T)\le C\,,\label{boundvar0}
\\
& \|\mu_t\|_*\le C_*+C \quad \hbox{for every }t\in[0,T] \,,\label{boundnorm0}
\end{eqnarray}
such that
\begin{equation}\label{convergence}
\mu^k_{t_1\dots t_m} \wto\mu_{t_1\dots t_m} \quad\hbox{weakly}^*\hbox{ in }GY(X 
;\Xi^m) 
\end{equation}
for every finite sequence $t_1,\dots,t_m$ in $\Theta$ with $0\le t_1 <\dots < t_m\le T$. 
\end{theorem}

\par
\begin{proof}
The proof is divided in several steps.

\par
\medskip
\noindent
{\it Step 1. Boundedness of $\mu^k_{t_1\dots t_m}$.}
We begin by proving that $\|\mu^k_t\|_*$ is bounded uniformly with respect to 
$t\in[0,T]$ and $k$. Let us fix $t<t_0$. By the compatibility condition (\ref{compatib})
\begin{eqnarray*}
& \langle |\xi|,\mu^k_t(x,\xi,\eta)\rangle -\langle|\xi_0|,\mu^k_{t_0}(x,\xi_0,\eta)\rangle = 
\langle|\xi|, \mu^k_{tt_0}(x,\xi,\xi_0,\eta)\rangle - \langle |\xi_0|, 
\mu^k_{tt_0}(x,\xi,\xi_0,\eta)\rangle \le \\
& \le \langle |\xi-\xi_0|, \mu_{tt_0}^k(x,\xi,\xi_0,\eta)\rangle \le
{\rm Var}(\mu^k;t,t_0)\le C\,.
\end{eqnarray*}
Thanks to Remark~\ref{M+}, from (\ref{boundvar}) and (\ref{boundnorm}) we obtain 
that
\begin{equation}\label{boundnorm2}
\sup_{k}\|\mu^k_t\|_*\le C_*+C
\end{equation}
for every $t\in[0,t_0)$.
A similar argument proves (\ref{boundnorm2}) when $t\in[t_0,T]$.

\par
By Lemma~\ref{900} and  (\ref{boundnorm2}) we obtain
\begin{equation}\label{boundnorm3}
\|\mu^k_{t_1\dots t_m}\|_*\le m(C_*+C)
\end{equation}
for every finite sequence $t_1,\dots,t_m$ with $t_1 <\dots < t_m$.

\par
\medskip
\noindent
{\it Step 2. Choice of the subsequence.}
Let $D$ be a countable dense subset of $[0,T]$ containing $0$. By the compactness
Theorem~\ref{thm:cpt}, using  (\ref{boundnorm3}) and a diagonal argument, we can
extract a subsequence, still denoted $\mu^k$, such that, for every $s_1,\dots, s_m$
in $D$ with $0\le s_1 <\dots < s_m\le T$, the sequence $\mu^k_{s_1\dots s_m}$
converges weakly$^*$ in $GY(X ;\Xi^m)$.

\par
\medskip
\noindent
{\it Step 3. Choice of $\Theta$.}
Let $V^k(t):={\rm Var}(\mu^k;0,t)$. Since $V^k$ is nondecreasing, by (\ref{boundvar}) 
and by Helly's Theorem there exists a subsequence, still denoted $V^k$, such that, for 
every $t\in[0,T]$, $V^k(t)\to V(t)$, where $V$ is a 
nondecreasing function on $[0,T]$ with values in $[0,C]$. Let
\begin{equation}\label{Theta}
\Theta:=\{0\}\cup\{t\in (0,T]:  \lim_{s\to t-}V(s)=V(t)\}\,.
\end{equation}

\par
\medskip
\noindent
{\it Step 4. Convergence and left continuity on $\Theta$.}
Let us fix two finite sequences $t_1,\dots,t_m$ and $s_1,\dots,s_m$ in $[0,T]$ such that
$0\le s_1 < t_1 <\dots <s_m < t_m \le T$. We want to estimate the difference 
$\mu^k_{t_1\dots t_m}-\mu^k_{s_1\dots s_m}$. Let $f\in 
C^\hom_L(X{\times}\Xi^m{\times}\R)$. Then there exists a constant $a$ such that
\begin{equation}\label{Lipschitz}
|f(x,\xi_{t_1},\dots,\xi_{t_m},\eta)-f(x,\xi_{s_1},\dots,\xi_{s_m},\eta)|\le a \sum_{i=1}^m 
|\xi_{t_i}-\xi_{s_i}|\,.
\end{equation}
By the compatibility condition (\ref{compatib}) we have the estimate
\begin{eqnarray*}
& \displaystyle\vphantom{\sum_i} |\langle f, \mu^k_{t_1\dots t_m} \rangle - \langle f, 
\mu^k_{s_1\dots s_m} \rangle | = 
\\
&=\! |\langle f(x,\xi_{t_1},\dots,\xi_{t_m},\eta)-
f(x,\xi_{s_1},\dots,\xi_{s_m},\eta), 
\mu^k_{s_1t_1\dots s_m t_m}(x,\xi_{s_1},\xi_{t_1},\dots,\xi_{s_m},\xi_{t_m},\eta) 
\rangle | \!\le
\\
& \displaystyle \le a \sum_{i=1}^m \langle | \xi_{t_i}-\xi_{s_i}|, \mu^k_{s_1t_1\dots 
s_m t_m}(x,\xi_{s_1},\xi_{t_1},\dots,\xi_{s_m},\xi_{t_m},\eta)\rangle =
\\
&  \displaystyle  = a \sum_{i=1}^m \langle | \xi_{t_i}-\xi_{s_i} |, 
\mu^k_{s_it_i}(x,\xi_{s_i},\xi_{t_i},\eta) \rangle \,,
\end{eqnarray*} 
which by (\ref{additive}) gives
\begin{equation}\label{estim}
|\langle f, \mu^k_{t_1\dots t_m} \rangle - \langle f, \mu^k_{s_1\dots s_m} \rangle |
\le a  \sum_{i=1}^m (V^k(t_i) -V^k(s_i))\,.
\end{equation} 
A simple modification of the proof shows that (\ref{estim}) holds even if
$0= s_1 = t_1< s_2 \le t_2 <\dots <s_m \le t_m \le T$. 

\par
If $t_1,\dots,t_m\in\Theta$ with $0\le t_1 <\dots < t_m \le T$, for every $\e$ we can
choose $s_1,\dots,s_m\in D$, with $0\le s_1 \le t_1< s_2 \le t_2 <\dots <s_m \le t_m \le 
T$,
such that $a\sum_i (V(t_i) -V(s_i)) <\e$. Using (\ref{estim}) we deduce that 
$|\langle f, \mu^k_{t_1\dots t_m} \rangle - \langle f, \mu^k_{s_1\dots s_m} \rangle | <\e$
for $k$ large enough. Since the sequence $\langle f, \mu^k_{s_1\dots s_m} \rangle$
converges, we have $|\langle f, \mu^k_{s_1\dots s_m} \rangle - 
\langle f, \mu^{k'}_{s_1\dots s_m} \rangle | <\e$ for $k, k'$ large enough.
It follows that $|\langle f, \mu^k_{t_1\dots t_m} \rangle - 
\langle f, \mu^{k'}_{t_1\dots t_m} \rangle | <3\e$ for $k, k'$ large enough,
hence $\langle f, \mu^k_{t_1\dots t_m} \rangle$ is a Cauchy sequence for every
$f\in C^\hom_L(X{\times}\Xi^m{\times}\R)$.
By (\ref{boundnorm3}) we deduce from Lemma~\ref{L-density} that, for every
$t_1,\dots,t_m$ in $\Theta$ with $0\le t_1 <\dots < t_m \le T$, the sequence
$\mu^k_{t_1\dots t_m}$ converges weakly$^*$ to some element $\mu_{t_1\dots t_m}$
of $GY(X ;\Xi^m)$ satisfying
\begin{equation}\label{boundnorm4}
\|\mu_{t_1\dots t_m}\|_*\le m(C_*+C)\,.
\end{equation}

\par
We observe that, given $t_1,\dots,t_m$ and $s_1\dots s_m$ in $\Theta$,
we can pass to the limit in (\ref{estim}) and obtain
\begin{equation}\label{estim2}
|\langle f, \mu_{t_1\dots t_m} \rangle - \langle f, \mu_{s_1\dots s_m} \rangle | 
\le a  \sum_{i=1}^m (V(t_i) -V(s_i))
\end{equation}
for every $f$ satisfying (\ref{Lipschitz}) and every pair of finite sequences
$t_1,\dots,t_m$ and $s_1,\dots,s_m$ in $\Theta$ such that
$s_1 \le t_1 <s_2\le t_2<\dots <s_m \le t_m$.
Using the definition \eqref{Theta} of $\Theta$ and Lemma~\ref{L-density}, we deduce 
from   
(\ref{boundnorm4}) and \eqref{estim2} that, for every $t_1,\dots,t_m$ in $\Theta$ with
$t_1 <\dots < t_m$, we have $\mu_{s_1\dots s_m}\wto \mu_{t_1\dots t_m}$
weakly$^*$ in $GY(X ;\Xi^m)$, as $s_i\to t_i$, $s_i\in\Theta$, and $s_i\le t_i$. 

\par
\medskip
\noindent
{\it Step 5. Extension to $[0,T]$.}
It remains to show that we can define $\mu_{t_1\dots t_m}$ when some $t_i$ does not
belong to $\Theta$, in such a way that the resulting system of generalized Young
measures satisfies the compatibility conditions, inequalities (\ref{boundvar0}) and
(\ref{boundnorm0}), and the continuity property (\ref{leftlim}).
To this purpose, it is enough to observe that, since $V$ has a finite limit from the left 
at each point, we have
\begin{equation}\label{V^-}
\lim_{k,k'\to\infty}\  \sum_{i=1}^m (V(s^k_i) -V(s^{k'}_i))=0
\end{equation}
for every sequence $(s^k_1,\dots,s^k_m)$ in $\Theta^m$ with 
$s^k_i\to t_i$, $s^k_i\le t_i$. Indeed, if $t_i\not\in\Theta$, we have
$$
V(s^k_i)\to V^-(t_i):=\lim_{\genfrac{}{}{0pt}2{\scriptstyle s\to t_i}{\scriptstyle s<t_i}}
V(s)\,.
$$
For these sequences $(s^k_1,\dots,s^k_m)$ we can deduce from estimate (\ref{estim2})
that $\langle f, \mu_{s^k_1\dots s^k_m} \rangle$ satisfies a Cauchy condition for every 
$f$ satisfying (\ref{Lipschitz}). By (\ref{boundnorm4}) we deduce from
Lemma~\ref{L-density} the existence of the weak$^*$-limit of $\mu_{s_1\dots s_m}$ as $s_i\to t_i$, 
$s_i\in\Theta$, and $s_i\le t_i$. We take such a weak$^*$ limit as the definition of
$\mu_{t_1\dots t_m}$. Clearly $\mu_{t_1\dots t_m}$ satisfies (\ref{boundnorm4}) and, 
by
construction, from \eqref{estim2}, we deduce that for every $f$ satisfying 
(\ref{Lipschitz})
and every pair of finite sequences $t_1,\dots,t_m$ and $s_1,\dots,s_m$ in $[0,T]$, with
$0\le s_1 \le t_1<s_2\le t_2< \dots <s_m \le t_m \le T$, there holds 
\begin{equation}\label{estim3}
|\langle f, \mu_{t_1\dots t_m} \rangle - \langle f, \mu_{s_1\dots s_m} \rangle | 
\le a \sum_{i=1}^m (V^-(t_i) -V^-(s_i))\,,
\end{equation}
where $V^-$ is the left-continuous representative of $V$ defined by \eqref{V^-}.
The continuity property (\ref{leftlim}) follows easily from \eqref{estim3} and from
Lemma~\ref{L-density}.

\par
For every finite sequence $t_1,\dots,t_m$ in $\Theta$ with $t_1 <\dots < t_m$
we have
$$
\sum_{i=1}^m \langle | \xi_i-\xi_{i-1} |, \mu^k_{t_{i-1}t_i}(x,\xi_{i-1},\xi_i,\eta) \rangle 
\le C\,.
$$
Passing to the limit as $k\to\infty$, we obtain 
$$
\sum_{i=1}^m \langle | \xi_i-\xi_{i-1} |, \mu_{t_{i-1}t_i}(x,\xi_{i-1},\xi_i,\eta) \rangle  
\le C
$$
whenever $t_1,\dots,t_m \in\Theta$. This restriction can be removed by an approximation
argument, and this proves (\ref{boundvar0}).

\par
The compatibility condition (\ref{compatib}) for $\mu^k$ implies that
$$
\langle f, \mu^k_{s_1\dots s_n}\rangle = \langle f\circ\pi^{t_1\dots t_m}_{s_1\dots s_n},
\mu^k_{t_1\dots t_m}\rangle
$$
for every $f\in C^\hom(X{\times}\Xi^h{\times}\R)$ and
every pair of finite sequences $s_1,\dots,s_n$ and $t_1,\dots,t_m$ in $[0,T]$ with 
$s_1<\dots<s_n$, $t_1<\dots<t_m$, and $\{s_1,\dots,s_n\}\subset \{t_1,\dots, t_m\}$. 
Passing to the limit as $k\to\infty$, we obtain 
$$
\langle f, \mu_{s_1\dots s_n}\rangle = \langle f\circ\pi^{t_1\dots t_m}_{s_1\dots s_h},
\mu_{t_1\dots t_m}\rangle\,,
$$
whenever $s_i$ and $t_j$ belong to $\Theta$. This restriction can be removed by an 
approximation argument, therefore $\mu\in SGY([0,T],X;\Xi)$.
\end{proof}

\par
We conclude this section by proving the lower semicontinuity of the $h$-variation.

\par
\begin{theorem}\label{semicontinuity} 
Let $T>0$ and let $\mu^k$ be a sequence in $SGY([0,T],X ;\Xi)$. 
Suppose that there exist a dense set $D\subset [0,T]$, and a left continuous 
$\mu\in SGY([0,T],X ;\Xi)$ such that
$$
\mu^k_{t_1\dots t_m} \wto\mu_{t_1\dots t_m} \quad\hbox{weakly}^*\hbox{ in }GY(X 
;\Xi^m) 
$$
for every finite sequence $t_1,\dots,t_m$ in $D$ with $t_1 <\dots < t_m$. Then
$$
 {\rm Var}_h(\mu;0,T)\le \liminf_{k\to\infty} {\rm Var}_h(\mu^k;0,T)
$$
for every positively one-homogeneous function $h\colon\Xi\to[0,+\infty)$ satisfying the 
triangle inequality.
\end{theorem}

\par
\begin{proof}
Let us fix $h$. For every finite sequence $t_1,\dots,t_m$ in $D$ with $t_1 <\dots < t_m$
we have
$$
\sum_{i=1}^m \langle h( \xi_i-\xi_{i-1} ), \mu^k_{t_{i-1}t_i}(x,\xi_{i-1},\xi_i,\eta) 
\rangle \le 
{\rm Var}_h(\mu^k;0,T)\,.
$$
Since $h$ is continuous (Remark~\ref{rem-Varh}), passing to the limit as $k\to\infty$ 
we obtain 
$$
\sum_{i=1}^m \langle h( \xi_i-\xi_{i-1} ), \mu_{t_{i-1}t_i}(x,\xi_{i-1},\xi_i,\eta) \rangle  
\le 
\liminf_{k\to\infty} {\rm Var}_h(\mu^k;0,T)
$$
whenever $t_1,\dots,t_m \in D$. The same inequality can be proved when 
$t_1,\dots,t_m \in [0,T]$ by an approximation argument, thanks to left continuity. The conclusion is 
obtained by taking the supremum with respect to $t_1,\dots,t_m$.
\end{proof}

\end{section}

\begin{section}{Weak$^*$ derivatives of systems with bounded variation}

In this section we introduce the notion of weak$^*$ derivative of a compatible system of 
generalized Young measures on the time interval $[0,T]$, with $T>0$, and prove that, if 
${\rm Var}(\mu;0,T)<+\infty$, then the weak$^*$ derivative exists at almost every 
$t\in [0,T]$.

\par
\begin{definition}\label{defdiffquot}
Given $\mu\in SGY([0,T],X;\Xi)$, the {\em difference quotient\/} of  $\mu$ between 
times $t_1$ and $t_2$, with $0\le t_1<t_2\le T$, is the element of $GY(X;\Xi)$ defined 
as the image
$$
q_{t_1t_2}(\mu_{t_1t_2})
$$
of $\mu_{t_1t_2}$ under the map $q_{t_1t_2}\colon 
X{\times}\Xi{\times}\Xi{\times}\R\to X{\times}\Xi{\times}\R$ defined by
$$
\textstyle
q_{t_1t_2}(x,\xi_1,\xi_2,\eta)=(x,\frac{\xi_2-\xi_1}{t_2-t_1},\eta)\,.
$$
\end{definition}

\par
\begin{remark}\label{rem673}
It follows from Definition~\ref{def:image3} that the difference quotient is characterized 
by the equality
$$
\langle f(x,\xi,\eta),q_{t_1t_2}(\mu_{t_1t_2})(x,\xi,\eta)\rangle=
\langle f(x,{\textstyle\frac{\xi_2-\xi_1}{t_2-
t_1}},\eta),\mu_{t_1t_2}(x,\xi_1,\xi_2,\eta)\rangle
$$
for every $f\in C^\hom(X{\times}\Xi{\times}\R)$.
\end{remark}

\par
\begin{remark}\label{rem674}
It follows from the definition of barycentre that
\begin{equation}\label{bardiffquot}
\bary(q_{t_1t_2}(\mu_{t_1t_2}))=\frac{\bary(\mu_{t_2})-\bary(\mu_{t_1})}{t_2-t_1}\,.
\end{equation}
In particular, if $\mu$ is the compatible system of generalized Young measures associated 
with a function $t\mapsto p(t)$ from $[0,T]$ into $M_b(X;\Xi)$ according to 
\eqref{delta1m}, then
\begin{equation}\label{pdiffq}
q_{t_1t_2}(\mu_{t_1t_2})=q_{t_1t_2}(\delta_{(p(t_1),p(t_2))})=\delta_{\frac{p(t_2)-
p(t_1)}{t_2-t_1}}
\end{equation}
(see Definition~\ref{GYMp}).
\end{remark}

\par
\begin{definition}\label{weak*der}
We say that $\mu\in SGY([0,T],X;\Xi)$ has a {\em weak$^*$ derivative\/} 
$\dot\mu_{t_0}$ at time $t_0\in[0,T]$ if $q_{tt_0}(\mu_{tt_0})\wto\dot\mu_{t_0}$ 
weakly$^*$ in  $GY(X;\Xi)$ as $t\to t_0^-$ and 
$q_{t_0t}(\mu_{t_0t})\wto\dot\mu_{t_0}$ weakly$^*$ in  $GY(X;\Xi)$ as $t\to 
t_0^+$, which is equivalent to
\begin{equation}\label{defder}
\begin{array}{c}
\displaystyle
\langle f(x,\xi_0,\eta),\dot\mu_{t_0}(x,\xi_0,\eta)\rangle=
\lim_{t\to t_0^-} 
\langle f(x,{\textstyle\frac{\xi-\xi_0}{t-t_0}},\eta),\mu_{tt_0}(x,\xi,\xi_0,\eta)\rangle=
\\
\displaystyle
=\lim_{t\to t_0^+} 
\langle f(x,{\textstyle\frac{\xi-\xi_0}{t-t_0}},\eta),\mu_{t_0t}(x,\xi_0,\xi,\eta)\rangle
\end{array}
\end{equation}
for every $f\in C^\hom(X{\times}\Xi{\times}\R)$.
\end{definition}

\par
\begin{remark}\label{weak*der2}
It follows from (\ref{pdiffq}) that, if $\mu$ is the compatible system of generalized 
Young measures associated with a function $t\mapsto p(t)$ from $[0,T]$ into 
$M_b(X;\Xi)$ according to \eqref{delta1m} and
$$
\frac{p(t)-p(t_0)}{t-t_0}\to \dot p(t_0)
$$
strongly in $M_b(X;\Xi)$ as $t\to t_0$, then $\mu$ has a weak$^*$ derivative at $t_0$ 
and
$$
\dot\mu_{t_0}=\delta_{\dot p(t_0)}\,.
$$
This is not true if
\begin{equation}\label{derpt}
\frac{p(t)-p(t_0)}{t-t_0}\wto \dot p(t_0)
\end{equation}
only in the weak$^*$ topology of $M_b(X;\Xi)$. However, using 
Remark~\ref{rem674}, in this case we obtain
$$
\bary(\dot\mu_{t_0})= \dot p(t_0)\,,
$$
if the weak$^*$ derivative of $\mu_{t_1\dots t_m}:=\delta_{(p(t_1),\dots,p(t_m))}$ 
exists at~$t_0$.

\par
An example where (\ref{derpt}) holds but $\dot\mu_{t_0}\neq \delta_{\dot p(t_0)}$, can 
be constructed in the following way. Let $T=2$, $X=[-1,1]$, $\Xi=\R$, let $\lambda$ be 
the Lebesgue measure, let $w\colon\R\to\R$ be the $2$-periodic function defined by
$$
w(x):=\begin{cases}
1 & \text{if } 2k\le x < 2k+1\,\text{ for some }\,k\in \Z\,,
\\
-1 & \text{if } 2k-1\le x<2k\,\text{ for some }\,k\in \Z\,.
\end{cases}
$$
For every $t\in[0,2]$ let $u(t)\in L^1(X)$ be the function defined by
$$
u(t,x):=\begin{cases}
(t-1)w(\frac{x}{t-1}) & \text{if } t\neq 1\,,
\\
0 & \text{if } t=1\,.
\end{cases}
$$
As $t\to 1$ we have
$$
\frac{u(t)-u(1)}{t-1} \wto 0 \qquad \hbox{weakly}^* \hbox{ in } M_b(X;\Xi)\,,
$$
while
$$
\delta_{\frac{u(t)-u(1)}{t-1}} \wto \tfrac12\delta_1 +\tfrac12 \delta_{-1}  \qquad 
\hbox{weakly}^* \hbox{ in } GY(X;\Xi)\,,
$$
which implies $\dot\mu_1=\tfrac12\delta_1 +\tfrac12 \delta_{-1}$.
\end{remark}

\par
\begin{remark}\label{weak*der3}
If $\mu\in SGY([0,T],X;\Xi)$ has a weak$^*$ derivative $\dot\mu_{t_0}$ at time 
$t_0\in[0,T]$, then
$$
\frac{\bary(\mu_t)-\bary(\mu_{t_0})}{t-t_0}\wto \bary(\dot\mu_{t_0})
$$
weakly$^*$ in $M_b(X;\Xi)$ as $t\to t_0$.
This follows from \eqref{bardiffquot} and Remark~\ref{bar11}.
\end{remark}

\par
The following theorem is the main result of this section.

\par
\begin{theorem}\label{BV-weakder}
Let $T>0$ and let $\mu\in SGY([0,T],X;\Xi)$ with ${\rm Var}(\mu;0,T)<+\infty$. Then 
the weak$^*$ derivative $\dot\mu_t$ exists for a.e.\ $t\in[0,T]$. Moreover, for every 
$f\in C^\hom(X{\times}\Xi{\times}\R)$ the function $t\mapsto \langle 
f,\dot\mu_t\rangle$ is integrable on $[0,T]$. Finally, if $h\colon\Xi\to[0,+\infty)$ is a 
positively one-homogeneous function satisfying the triangle inequality, then 
\begin{equation}\label{BV-Varder}
\int_a^b\langle h(\xi),\dot\mu_t(x,\xi,\eta)\rangle\,dt \le {\rm Var}_h(\mu;a,b)
\end{equation}
for every $a,\,b\in[0,T]$ with $a\le b$.
\end{theorem}

\par
\begin{proof}
Some ideas of this proof are borrowed from the proof of 
\cite[Theorem~4.1.1]{Amb-Til} 
on the existence of the metric derivative of a Lipschitz curve.

\par
\medskip
\noindent
{\it Step 1. Boundedness of the difference quotients.\/}
By Remark~\ref{M+} and by \eqref{suppmu} for every $t_1,\,t_2\in[0,T]$, with 
$t_1<t_2$, we have 
$$
\textstyle
\|q_{t_1t_2}(\mu_{t_1t_2})\|_*\le \frac{1}{t_2-t_1} \langle |\xi_2-
\xi_1|,\mu_{t_1t_2}(x,\xi_1,\xi_2,\eta)\rangle + 
 \langle \eta,\mu_{t_1t_2}(x,\xi_1,\xi_2,\eta)\rangle\,.
$$
Let $V\colon [0,T]\to[0,+\infty)$ be the nondecreasing function defined by
\begin{equation}\label{BV-V(t)}
V(t):={\rm Var}(\mu;0,t)\,.
\end{equation}
By (\ref{muproj}) and (\ref{additive}) we conclude that
$$
\textstyle \|q_{t_1t_2}(\mu_{t_1t_2})\|_*\le  \frac{V(t_2)-V(t_1)}{t_2-t_1}+ 
\lambda(X)\,.
$$
Let $t_0\in[0,T]$ be a point where the derivative of $V$ exists. By the previous inequality 
we have that 
$$
\textstyle \|q_{tt_0}(\mu_{tt_0})\|_* \qquad \hbox{and} \qquad \textstyle 
\|q_{t_0t}(\mu_{t_0t})\|_*
$$
are bounded uniformly with respect to $t$. 
By the separability of $C^\hom(X{\times}\Xi{\times}\R)$ there exists a countable dense 
subset ${\mathcal F}$ of the set $C^\hom_\triangle(X{\times}\Xi{\times}\R)$ 
introduced in Definition~\ref{ChomT}. Therefore, since ${\mathcal F}$ is dense in 
$C^\hom(X{\times}\Xi{\times}\R)$ (see Lemma~\ref{convhom}), to prove the existence 
of the weak$^*$ derivative of $\mu$ at $t_0$ it is enough to show that
\begin{equation}\label{BV-der+-}
\lim_{t\to t_0^-} 
\langle f(x,{\textstyle\frac{\xi-\xi_0}{t-t_0}},\eta),\mu_{tt_0}(x,\xi,\xi_0,\eta)\rangle=
\lim_{t\to t_0^+} 
\langle f(x,{\textstyle\frac{\xi-\xi_0}{t-t_0}},\eta),\mu_{t_0t}(x,\xi_0,\xi,\eta)\rangle
\end{equation}
for every $f\in {\mathcal F}$.

\par
\medskip
\noindent
{\it Step 2. Some auxiliary functions.\/}
In order to prove (\ref{BV-der+-}), let us fix $f\in {\mathcal F}$
and let $\tau_i$ be a countable dense sequence in $[0,T]$. {}For every $i$ we define
\begin{equation}\label{BV-phii}
\varphi_i^f(t):=
\begin{cases}
\langle f(x, \xi-\zeta_i, (t-\tau_i)\eta), \mu_{t\tau_i}(x,\xi,\zeta_i,\eta)\rangle& \text{if 
}t<\tau_i\,,
\\
0& \text{if }t=\tau_i\,,
\\
\langle f(x, \xi-\zeta_i, (t-\tau_i)\eta), \mu_{\tau_it}(x,\zeta_i,\xi,\eta)\rangle& \text{if 
}t>\tau_i\,.
\end{cases}
\end{equation}

\par
Let us prove that $\varphi_i^f$ has bounded variation. Let us fix $t_1,\,t_2\in[0,T]$, with 
$t_1<t_2$. We consider first the case $t_1<\tau_i<t_2$. By the compatibility condition 
(\ref{compatib}) we have
$$
|\varphi_i^f(t_2)-\varphi_i^f(t_1)|\le 
\langle |f(x,\xi_2-\zeta_i,(t_2-\tau_i)\eta)-f(x,\xi_1-\zeta_i,(t_1-\tau_i)\eta)|, 
\mu_{t_1\tau_it_2}(x,\zeta_i,\xi_1,\xi_2,\eta)\rangle\,.
$$
Since, by Remark~\ref{rem:tri+lip}, 
$$
|f(x,\xi_2-\zeta_i,(t_2-\tau_i)\eta)-f(x,\xi_1-\zeta_i,(t_1-\tau_i)\eta)|\le
(|\xi_2-\xi_1|+(t_2-t_1)|\eta|)\,\|f\|_\hom\,,
$$
using again (\ref{compatib}) we obtain
$$
|\varphi_i^f(t_2)-\varphi_i^f(t_1)|\le 
\big(\langle |\xi_2-\xi_1|,\mu_{t_1t_2}(x,\xi_1,\xi_2,\eta)\rangle +  
(t_2-t_1)\,\langle |\eta| ,\mu_{t_1t_2}(x,\xi_1,\xi_2,\eta)\rangle\big)\,\|f\|_\hom \,.
$$
The same inequality can be proved when $\tau_i\le t_1$ or $\tau_i\ge t_2$.
By \eqref{suppmu}, (\ref{muproj}),  (\ref{additive}), and (\ref{BV-V(t)}) we conclude 
that
\begin{equation}\label{BV-lip30}
|\varphi_i^f(t_2)-\varphi_i^f(t_1)| \le (V(t_2)-V(t_1)+(t_2-t_1)\lambda(X))\, \|f\|_\hom\,.
\end{equation}

\par
We now prove that for every $t_1,\,t_2\in[0,T]$, with $t_1<t_2$, we have
\begin{equation}\label{BV-estbelow60}
\varphi_i^f(t_2)-\varphi_i^f(t_1)\le 
\langle f(x,\xi_2-\xi_1, (t_2-t_1)\eta), \mu_{t_1t_2}(x,\xi_1,\xi_2,\eta)\rangle\,.
\end{equation}
We consider first the case $t_1<\tau_i<t_2$. By (\ref{compatib}) and (\ref{BV-phii}) we 
have
$$
\varphi_i^f(t_2)=\langle f(x,\xi_2-\xi_1+\xi_1-\zeta_i, (t_2-t_1+t_1-\tau_i)\eta),
\mu_{t_1\tau_it_2}(x,\xi_1,\zeta_i,\xi_2,\eta)\rangle\,.
$$
{}From the triangle inequality and from (\ref{compatib}) we get
\begin{equation}\label{BV-estbelow61}
\begin{array}{c}
\varphi_i^f(t_2)\le \langle f(x,\xi_2-\xi_1, (t_2-t_1)\eta),
\mu_{t_1t_2}(x,\xi_1,\xi_2,\eta)\rangle+{}
\\
{}+\langle f(x,\xi_1-\zeta_i, (t_1-\tau_i)\eta),
\mu_{t_1\tau_i}(x,\xi_1,\zeta_i,\eta)\rangle\,,
\end{array}
\end{equation}
which gives \eqref{BV-estbelow60} by (\ref{BV-phii}). The proof in the cases 
$\tau_i\le t_1$ and $\tau_i\ge t_2$ is similar.

\par
Let $W\colon[0,T]\to \R$ be the increasing function defined by
\begin{equation}\label{W(t)}
W(t):=V(t)+t\,\lambda(X)
\end{equation} 
and let $\sigma\colon[0,W(T)]\to[0,T]$ be the nondecreasing function defined by
$$
\sigma(s):=\inf\{t\in[0,T]:W(t)\ge s\}\,.
$$
It is easy to see that 
\begin{equation}\label{sigmaW}
\sigma(W(t))=t \qquad \hbox{for every } t\in[0,T]\,.
\end{equation} 
As $W(t_2)-W(t_1)\ge (t_2-t_1)\lambda(X)$ for every $t_1<t_2$, we have 
\begin{equation}\label{Lip-sigma}
0\le \sigma(s_2)-\sigma(s_1)\le (s_2-s_1)/\lambda(X)
\end{equation} 
for every $s_1<s_2$, hence $\sigma$ is Lipschitz continuous. 

\par
By \eqref{BV-lip30} and \eqref{sigmaW} we have
$$
| (\varphi_i^f\circ\sigma)(s_2)-  (\varphi_i^f\circ\sigma)(s_1) |\le |s_2-s_1|\, \|f\|_\hom
$$
for every $s_1,s_2\in W([0,T])$. Therefore, there exists a function $\psi_i^f
\colon [0,W(T)]\to\R$ such that $\psi_i^f(s)=(\varphi_i^f\circ \sigma)(s)$ for every $s\in 
W([0,T])$ and 
\begin{equation}\label{Lip-psi}
| \psi_i^f(s_2)-  \psi_i^f
(s_1) |\le |s_2-s_1|\, \|f\|_\hom
\end{equation}
for every $s_1,s_2\in [0,W(T)]$. For every $s_0\in [0,W(T)]$ let 
$$
\dot\psi_i^f(s_0)=\limsup_{s\to s_0} \frac{\psi_i^f
(s)-\psi_i^f
(s_0)}{s-s_0}\,.
$$
By \eqref{Lip-psi} we have $|\dot\psi_i^f(s_0)|\le \|f\|_\hom$, and by 
Lebesgue's Differentiation Theorem the limsup is a limit for a.e.\ $s_0\in[0,W(T)]$.
Finally, let $\omega^f\colon[0,W(T)]\to \R$ be the function defined by
\begin{equation}\label{BV-psiphi'}
\omega^f(s):=\sup_i \dot\psi_i^f(s)\,.
\end{equation}
By the bound on $\dot\psi_i^f$ we have
\begin{equation}\label{378}
|\omega^f(s)|\le \|f\|_\hom 
\end{equation}
for every $s\in [0,W(T)]$.

\par
\medskip
\noindent
{\it Step 3. The exceptional set.\/}
Let $\Lone$ be the Lebesgue measure on $\R$.
By \eqref{Lip-sigma} and \eqref{Lip-psi} there exists a measurable set 
$N\subset [0,W(T)]$, with $\Lone(N)=0$, such that each point of 
$[0,W(T)]\setminus N$ is a Lebesgue point of $\omega^f$ for every 
$f\in{\mathcal F}$ and a differentiability point for $\psi_i^f$
for every $f\in{\mathcal F}$ and for every~$i$. Let $N_W$ be the set of points of 
$[0,T]$ where the derivative $\dot W$ of $W$ does not exist. By Lebesgue's
Differentiation Theorem we have $\Lone(N_W)=0$. Since $\sigma$ is 
Lipschitz continuous and $W^{-1}(N)=\sigma(N\cap W([0,T]))$ by \eqref{sigmaW}, 
we have that $\Lone(W^{-1}(N))=0$, hence
\begin{equation}\label{Lone0}
\Lone(N_W\cup W^{-1}(N))=0\,.
\end{equation}

\par
\medskip
\noindent
{\it Step 4. Estimate from below.\/}
Let us fix $t_0\not\in N_W\cup W^{-1}(N)$, with $0<t_0<T$, and let $s_0=W(t_0)$. 
As $\varphi_i^f(t)=\psi_i^f
(W(t))$, from (\ref{BV-estbelow60}) we obtain
$$
\psi_i^f
(W(t_2))-\psi_i^f(W(t_1))\le 
\langle f(x,\xi_2-\xi_1, (t_2-t_1)\eta), \mu_{t_1t_2}(x,\xi_1,\xi_2,\eta)\rangle
$$
for every $t_1,t_2\in [0,T]$ with $t_1<t_2$. This implies
\begin{eqnarray*}
&\displaystyle
\dot\psi_i^f(W(t_0))\, \dot W(t_0)\le \liminf_{t\to t_0^-}  \, 
\langle f(x,{\textstyle\frac{\xi-\xi_0}{t-t_0}},\eta),\mu_{tt_0}(x,\xi,\xi_0,\eta)\rangle\,,
\\
&\displaystyle
\dot\psi_i^f(W(t_0))\, \dot W(t_0) \le 
\liminf_{t\to t_0^+} \, 
\langle f(x,{\textstyle\frac{\xi-\xi_0}{t-t_0}},\eta),\mu_{t_0t}(x,\xi_0,\xi,\eta)\rangle
\end{eqnarray*}
for every $i$, which by \eqref{BV-psiphi'} gives
\begin{eqnarray}
&\displaystyle
\omega^f(W(t_0))\, \dot W(t_0)\le \liminf_{t\to t_0^-}  \, 
\langle f(x,{\textstyle\frac{\xi-\xi_0}{t-t_0}},\eta),\mu_{tt_0}(x,\xi,\xi_0,\eta)\rangle\,,
\label{BV-estbelow62}
\\
&\displaystyle
\omega^f(W(t_0))\, \dot W(t_0) \le 
\liminf_{t\to t_0^+} \, 
\langle f(x,{\textstyle\frac{\xi-\xi_0}{t-t_0}},\eta),\mu_{t_0t}(x,\xi_0,\xi,\eta)\rangle\,.
\label{BV-estbelow63}
\end{eqnarray}

\par
\medskip
\noindent
{\it Step 5. Estimate from above.\/}
To prove the opposite inequality we show that
\begin{equation}\label{BV-estabove60}
\langle f(x,\xi_2-\xi_1, (t_2-t_1)\eta), \mu_{t_1t_2}(x,\xi_1,\xi_2,\eta)\rangle
\le \int_{W(t_1)}^{W(t_2)} \omega^f(s)\,ds
\end{equation}
for every $t_1,t_2$, with $0<t_1<t_2<T$, such that $W$ is continuous at $t_1$ or 
$t_2$. We prove \eqref{BV-estabove60} only when $W$ is continuous at $t_1$, the 
other case being analogous.
{}For every $\e>0$ there exists $i$ such that $\tau_i<t_1$ and 
\begin{equation}\label{cont-W}
W(t_1)-W(\tau_i)<\e\,.
\end{equation} 
As $\psi_i^f$ is Lipschitz, from the compatibility condition (\ref{compatib}) we get
\begin{equation}\label{BV-estabove64}
\begin{array}{c}
\displaystyle
\int_{W(t_1)}^{W(t_2)} \omega^f(s)\,ds \ge \int_{W(t_1)}^{W(t_2)} \dot\psi_i^f(s)\,ds 
=
\psi_i^f(W(t_2))-\psi_i^f(W(t_1))=
\\
\displaystyle\vphantom{\int_{t_1}^{t_2}}
 =
\varphi_i^f(t_2)-\varphi_i^f(t_1)
=\langle f(x,\xi_2-\zeta_i,(t_2-\tau_i)\eta), 
\mu_{\tau_it_1t_2}(x,\zeta_i,\xi_1,\xi_2,\eta)\rangle-{}
\\
\displaystyle\vphantom{\int_{t_1}^{t_2}}
{}-\langle f(x,\xi_1-\zeta_i,(t_1-\tau_i)\eta), 
\mu_{\tau_it_1t_2}(x,\zeta_i,\xi_1,\xi_2,\eta)\rangle\,.
\end{array}
\end{equation}
Using Remark~\ref{rem:tri+lip} we obtain
$$
\begin{array}{c}
f(x,\xi_2-\zeta_i,(t_2-\tau_i)\eta)\ge 
f(x,\xi_2-\xi_1,(t_2-t_1)\eta) - (|\xi_1-\zeta_i|+ (t_1-\tau_i)|\eta|)\,\|f \|_\hom\,,
\\
-f(x,\xi_1-\zeta_i,(t_1-\tau_i)\eta)\ge 
 - (|\xi_1-\zeta_i|+ (t_1-\tau_i)|\eta|)\,\|f \|_\hom\,,
\end{array}
$$
so that, using again (\ref{muproj}) and (\ref{compatib}), inequality
(\ref{BV-estabove64})
and the definition of $W$ give
$$
\int_{W(t_1)}^{W(t_2)} \omega^f(s)\,ds \ge
\langle f(x,\xi_2-\xi_1,(t_2-t_1)\eta), \mu_{t_1t_2}(x,\xi_1,\xi_2,\eta)\rangle- 
2(W(t_1) - W(\tau_i)) \|f\|_\hom \,.
$$
By (\ref{cont-W}) we conclude that
$$
\int_{W(t_1)}^{W(t_2)} \omega^f(s)\,ds \ge
\langle f(x,\xi_2-\xi_1,(t_2-t_1)\eta), \mu_{t_1t_2}(x,\xi_1,\xi_2,\eta)\rangle-
2\e\| f\|_\hom \,.
$$
As $\e>0$ is arbitrary, this proves~(\ref{BV-estabove60}).

\par
Since $W$ is differentiable at $t_0$ and $W(t_0)$ is a Lebesgue point of $\omega^f$, 
inequality (\ref{BV-estabove60}) implies
\begin{eqnarray*}
&\displaystyle
\limsup_{t\to t_0^-}  \, 
\langle f(x,{\textstyle\frac{\xi-\xi_0}{t-t_0}},\eta),\mu_{tt_0}(x,\xi,\xi_0,\eta)\rangle\le 
\omega^f(W(t_0))\, \dot W(t_0) \,,
\\
&\displaystyle
\limsup_{t\to t_0^+} \, 
\langle f(x,{\textstyle\frac{\xi-\xi_0}{t-t_0}},\eta),\mu_{t_0t}(x,\xi_0,\xi,\eta)\rangle\le 
\omega^f(W(t_0))\,\dot W(t_0) \,,
\end{eqnarray*}
which, together with (\ref{BV-estbelow62}) and (\ref{BV-estbelow63}), give
\begin{eqnarray*}
&\displaystyle
\lim_{t\to t_0^-} 
\langle f(x,{\textstyle\frac{\xi-\xi_0}{t-t_0}},\eta),\mu_{tt_0}(x,\xi,\xi_0,\eta)\rangle=
\omega^f(W(t_0))\, \dot W(t_0)\,,
\\
&\displaystyle
\lim_{t\to t_0^+} 
\langle f(x,{\textstyle\frac{\xi-\xi_0}{t-t_0}},\eta),\mu_{t_0t}(x,\xi_0,\xi,\eta)\rangle=
\omega^f(W(t_0))\, \dot W(t_0)\,.
\end{eqnarray*}
By (\ref{Lone0}) this proves (\ref{BV-der+-}) and concludes the proof of the existence 
of
the weak$^*$ derivative $\dot\mu_{t_0}$ for a.e.\ $t_0\in[0,T]$. Moreover it shows that
\begin{equation}\label{der+-2}
\langle f(x,\xi,\eta),\dot\mu_{t_0}(x,\xi,\eta)\rangle=\omega^f(W(t_0))\, \dot W(t_0)
\end{equation}
for every $f\in {\mathcal F}$ and for a.e.\ $t_0\in[0,T]$.

\par
\medskip
\noindent
{\it Step 6. Integrability of $t\mapsto\langle f,\dot\mu_t\rangle$.\/}
To prove the measurability of this function for every
$f\in C^\hom(X{\times}\Xi{\times}\R)$, we fix a sequence $\e_k$ of positive
numbers converging to $0$ and a function $f\in C^\hom_L(X{\times}\Xi{\times}\R)$.
By Lemma~\ref{lm:bv-mut} the function 
$$
t\mapsto \langle f(x,\tfrac{\xi'-\xi}{\e_k},\eta),\mu_{t,t+\e_k}(x,\xi,\xi',\eta)\rangle
$$ 
is measurable on $[0, T-\e_k]$. Since it converges to $t\mapsto \langle 
f,\dot\mu_t\rangle$ for a.e.\ $t\in[0,T]$, we conclude that this function is measurable on 
$[0,T]$. The same property can be proved for an arbitrary $f\in 
C^\hom(X{\times}\Xi{\times}\R)$ by approximation, thanks to Lemma~\ref{L-density}.

\par
By \eqref{378} and \eqref{der+-2} we have
$$
|\langle f,\dot\mu_t\rangle| \le \dot W(t)\,\|f\|_\hom
$$
for every $f\in {\mathcal F}$. The same inequality holds for
any $f\in C^\hom(X{\times}\Xi{\times}\R)$ by the density of ${\mathcal F}$ (see 
Lemma~\ref{convhom}). Since $ \dot W$ is integrable, this concludes the proof of the 
integrability of $t\mapsto\langle f,\dot\mu_t\rangle$ on $[0,T]$.

\par
\medskip
\noindent
{\it Step 7. Estimate for ${\rm Var}_h(\mu;a,b)$.\/}
Let $h\colon\Xi\to[0,+\infty)$ be a positively one-homogeneous function satisfying the 
triangle inequality. Since the function $t\mapsto{\rm Var}_h(\mu;a,t)$ is nondecreasing 
on $[a,b]$, by Lebesgue Differentiation Theorem it is differentiable for a.e.\ $t\in[a,b]$ 
and
\begin{equation}\label{BV-M011}
\int_a^b \frac{d}{dt}{\rm Var}_h(\mu;a,t)\, dt \le {\rm Var}_h(\mu;a,b)\,.
\end{equation}
Let $t_0\in(a,b)$ be a point where $t\mapsto{\rm Var}_h(\mu;a,t)$ is 
differentiable and the weak$^*$ derivative $\dot\mu_{t_0}$ exists.
By the definition of ${\rm Var}_h$ for every $t\in(t_0,b)$ we have
$$
{\rm Var}_h(\mu;a,t_0)+ \langle h(\xi-\xi_0), \mu_{t_0t}(x,\xi_0,\xi,\eta)\rangle \le
{\rm Var}_h(\mu;a,t)\,.
$$
Since $h$ is positively homogeneous of degree one, we obtain
$$
 \langle h({\textstyle\frac{\xi-\xi_0}{t-t_0}}), \mu_{t_0t}(x,\xi_0,\xi,\eta)\rangle
\le
\frac{{\rm Var}_h(\mu;a,t)-{\rm Var}_h(\mu;a,t_0)}{t-t_0}\,.
$$
{}From (\ref{defder}) we deduce that
$$
\langle h(\xi),\dot\mu_{t_0}(x,\xi,\eta)\rangle\le 
\frac{d}{dt}{\rm Var}_h(\mu;a,t)\Big|_{t=t_0}\,.
$$
Since this inequality holds for a.e.\ $t_0\in[a,b]$, from \eqref{BV-M011} we obtain
$$
 \int_a^b\langle h(\xi),\dot\mu_t(x,\xi,\eta)\rangle\,dt\le {\rm Var}_h(\mu;a,b)\,,
$$
which concludes the proof of~(\ref{BV-Varder}).
\end{proof}
\end{section}

\begin{section}{Absolute continuity}

In this section we introduce the notion of absolutely continuous system of generalized 
Young measures on the time interval $[0,T]$, with $T>0$, and prove that for these 
systems the $h$-variation can be computed using the weak$^*$ derivative by the formula 
\begin{equation}\label{Varhmut}
{\rm Var}_h(\mu;a,b)=\int_a^b\langle h(\xi),\dot\mu_t(x,\xi,\eta)\rangle\,dt
\end{equation}
for every $a,\,b\in[0,T]$, with $a<b$.

\par
\begin{definition}\label{abscont}
We say that a compatible system of generalized Young measures $\mu\in 
SGY([0,T],X;\Xi)$ is {\it absolutely continuous on\/} $[0,T]$ if for every $\e>0$ there 
exists $\delta>0$ such that
\begin{equation}\label{abscon3}
\sum_{i=1}^k\langle|\xi_2-\xi_1|,\mu_{a_ib_i}(x,\xi_1,\xi_2,\eta)\rangle\le\e
\end{equation}
for every finite family $(a_1, b_1),\dots,(a_k,b_k)$ of nonoverlapping open intervals in 
$[0,T]$ with
$$
\sum_{i=1}^k(b_i-a_i)\le\delta\,.
$$
\end{definition}

\par
\begin{remark}\label{abscont103}
It follows from Definition~\ref{GYMp} that, if $\mu$ is the compatible system of 
generalized Young measures associated with a function $t\mapsto p(t)$ from $[0,T]$ into 
$M_b(X;\Xi)$ according to \eqref{delta1m}, then $\mu$ is absolutely continuous on 
$[0,T]$ if and only if 
$t\mapsto p(t)$ is absolutely continuous on $[0,T]$ in the usual sense of functions with 
values in a Banach space.

\par
If $t\mapsto u(t)$ is an absolutely continuous function from $[0,T]$ into $L^r(X;\Xi)$ 
for some $r>1$, then the derivative $\dot u(t)$, defined as the strong $L^r$ limit of the 
difference quotients, exists at a.e.\ $t\in [0,T]$ (see, e.g., \cite[Appendix]{Bre}). By 
Remark~\ref{weak*der2} it follows that, if $\mu$ is the compatible system of generalized 
Young measures associated with  $t\mapsto u(t)$ according to \eqref{delta1m}, then 
$\dot \mu_t=\delta_{\dot u(t)}$ for a.e.\ $t\in [0,T]$, and \eqref{Varhmut} follows from the 
classical theory (see, e.g., \cite[Appendix]{Bre}).

\par
If $t\mapsto p(t)$ is an absolutely continuous function  with values in $M_b(X;\Xi)$, 
then the derivative $\dot p(t)$, defined as the weak$^*$ limit of the difference quotients, 
exists at a.e.\ $t\in [0,T]$ (see \cite[Appendix]{DM-DeS-Mor}). This is not enough to 
guarantee that $\dot \mu_t=\delta_{\dot p(t)}$ for a.e.\ $t\in [0,T]$ when $\mu$ is the 
compatible system of generalized Young measure associated with  $t\mapsto p(t)$ (see 
Remark~\ref{weak*der2}). Therefore, in this case \eqref{Varhmut} cannot be obtained 
directly from known results.
\end{remark}

\par
\begin{remark}\label{abscont2}
As in the classical case, it is easy to see that, if $\mu\in SGY([0,T],X;\Xi)$ is absolutely 
continuous on $[0,T]$, then ${\rm Var}(\mu;0,T)<+\infty$. In this case, if 
$V\colon[0,T]\to[0,+\infty)$ is the nondecreasing function defined by
$$
V(t):={\rm Var}(\mu;0,t)\,,
$$
then for every $\e>0$
$$
\sum_{i=1}^k(V(b_i)-V(a_i))\le\e
$$
for every finite family $(a_1, b_1),\dots,(a_k,b_k)$ of nonoverlapping open intervals in 
$[0,T]$ with
$$
\sum_{i=1}^k(b_i-a_i)\le\delta\,,
$$
where $\delta$ is the constant in the definition of the absolute continuity of $\mu$. In 
particular, $V$ is absolutely continuous on $[0,T]$.
\end{remark}

\par
\begin{theorem}\label{weakder}
Suppose that $\mu\in SGY([0,T],X;\Xi)$ is absolutely continuous on $[0,T]$ and that 
$h\colon\Xi\to[0,+\infty)$ is positively one-homogeneous and satisfies the triangle 
inequality. Then 
$$
{\rm Var}_h(\mu;a,b)=\int_a^b\langle h(\xi),\dot\mu_t(x,\xi,\eta)\rangle\,dt
$$
for every $a,\,b\in[0,T]$ with $a\le b$.
\end{theorem}

\par
\begin{proof} Let $W$ be defined by \eqref{W(t)}. By Remark~\ref{abscont2} $W$ is
absolutely continuous on $[0,T]$. By Remark~\ref{rem-Varh} the function
$f(x,\xi,\eta):=h(\xi)$ belongs to $C^\hom_\triangle(X{\times}\Xi{\times}\R)$.
Therefore, we can add this function to the set ${\mathcal F}$ introduced in Step~1 of
the proof of Theorem~\ref{BV-weakder} and we can consider the corresponding
function $\omega^h\colon[0,W(T)]\to \R$ defined by \eqref{BV-psiphi'}.
By \eqref{BV-estabove60} and (\ref{der+-2}) we have
\begin{eqnarray}
&\displaystyle
\langle h(\xi_2-\xi_1), \mu_{t_1t_2}(x,\xi_1,\xi_2,\eta)\rangle
\le \int_{W(t_1)}^{W(t_2)} \omega^h(s)\,ds\,,
\label{estaboveh}
\\
&
\langle h(\xi),\dot\mu_t(x,\xi,\eta)\rangle=\omega^h(W(t))\,\dot W(t)
\label{der+-3}
\end{eqnarray}
for every $t_1,\,t_2\in[0,T]$, with $t_1<t_2$, and for a.e.\ $t\in[0,T]$.

\par
By the definition of ${\rm Var}_h(\mu;a,b)$, inequality (\ref{estaboveh})  implies that
\begin{equation}\label{Varder<}
{\rm Var}_h(\mu;a,b)\le \int_{W(a)}^{W(b)}\omega^h(s)\,ds
\end{equation}
for every $a,\,b\in[0,T]$, with $a\le b$. On the other hand, since $W$ is absolutely 
continuous on $[0,T]$, we have
\begin{equation}\label{Varder<2}
\int_{W(a)}^{W(b)}\omega^h(s)\,ds=\int_a^b \omega^h(W(t))\,W'(t)\,dt=
\int_a^b\langle h(\xi),\dot\mu_t(x,\xi,\eta)\rangle\,dt\,,
\end{equation}
where the last equality follows from~(\ref{der+-3}). The conclusion follows now from 
\eqref{BV-Varder},
\eqref{Varder<}, and~\eqref{Varder<2}.
\end{proof}

\end{section}

\bigskip

\noindent {\bf Acknowledgments.} { This work is part of the project 
``Calculus of Variations" 2004, 
supported by the Italian Ministry of Education, University, and Research,
and of the project  ``Mathematical Challenges in Nanomechanics at the Interface between 
Atomistic and Continuum Models'' supported by INdAM.}

\bigskip
\bigskip

{
}

\end{document}